\documentclass[11pt, letterpaper, oneside]{amsart}

\usepackage{latexsym, amsmath, amssymb, amsfonts, amscd}
\usepackage{amsthm}
\usepackage{t1enc}
\usepackage[mathscr]{eucal}
\usepackage{indentfirst}
\usepackage{graphicx, pb-diagram}
\usepackage{fancyhdr}
\usepackage{fancybox}
\usepackage{enumerate}
\usepackage[all]{xy}
\usepackage{url}

\theoremstyle{plain}
\newtheorem{thm}{Theorem}[section]

\newtheorem{prop}[thm]{Proposition}
\newtheorem{lemma}[thm]{Lemma}

\theoremstyle{definition}
\newtheorem{defi}[thm]{Definition}

\theoremstyle{remark}
\newtheorem{remark}[thm]{Remark}
\newtheorem{ep}[thm]{Example}

\newcommand{\CC}{\ensuremath{\mathbb C}}
\newcommand{\RR}{\ensuremath{\mathbb R}}
\newcommand{\g}{\ensuremath{\frak{g}}}
\newcommand{\h}{\ensuremath{\frak{h}}}

\newcommand{\bt}{\mathbf{t}}                  
\newcommand{\bs}{\mathbf{s}}                  

\newcommand{\Gs}{\Gamma_s}

\newcommand{\cO}{\mathcal{O}}
\newcommand{\cF}{\mathcal{F}}
\newcommand{\cI}{\mathcal{I}}
\newcommand{\tP}{\tilde{P}}
\newcommand{\hP}{\hat{P}}

\newcommand{\sh}{\sharp}

\newcommand{\sctc}{\bs^{-1}(C)\cap \bt^{-1}(C)}
\newcommand{\stpttp}{\bs^{-1}(\tP)\cap \bt^{-1}(\tP)}

\begin{document}

\title{Coisotropic embeddings in Poisson manifolds}
\author{A.S. Cattaneo and M. Zambon}
\address{Institut f\"ur Mathematik, Universit\"at Z\"urich-Irchel, Winterthurerstr. 190, CH-8057 Z\"urich, Switzerland}
\email{alberto.cattaneo@math.unizh.ch, zambon@math.unizh.ch}

\begin{abstract}
We consider existence and uniqueness of two kinds of coisotropic
embeddings and deduce the existence of deformation quantizations
of certain Poisson algebras of basic functions. First we show that
any submanifold of a Poisson manifold  satisfying a certain
constant rank condition, already considered by Calvo and Falceto
\cite{CaFa1}, sits coisotropically inside some larger cosymplectic
submanifold, which is naturally endowed with a Poisson structure.
Then we give conditions under which a Dirac manifold can be
embedded coisotropically in a Poisson manifold, extending a
classical theorem of Gotay.
\end{abstract}

\maketitle
\tableofcontents

\section{Introduction}\label{intro}

The following two results in symplectic geometry are well known.
First: a submanifold $C$ of a symplectic manifold $(M,\Omega)$ is
contained coisotropically in some symplectic submanifold of $M$
iff the pullback of $\Omega$ to $C$ has constant rank; see Marle's work  \cite{Ma}. Second: a
manifold endowed with a closed 2-form $\omega$ can be embedded
coisotropically into a symplectic manifold $(M,\Omega)$ so that
$i^*\Omega=\omega$ (where $i$ is the embedding) iff $\omega$ has
constant rank; see Gotay's work \cite{Go}. 

In this paper we extend these results to
the setting of Poisson geometry. Recall that $P$ is a Poisson manifold if it is endowed with a bivector field $\Pi \in \Gamma(\wedge^2 TP)$ satisfying the Schouten-bracket condition $[\Pi,\Pi]=0$. A submanifold $C$ of $(P,\Pi)$ is coisotropic if $\sharp
N^*C\subset TC$, where the conormal bundle $N^*C$ is defined as the annihilator of $TC$ in $TP|_C$ and $\sharp\colon T^*P\rightarrow TP$ is the contraction with the
bivector $\Pi$. Coisotropic submanifolds appear naturally; for instance 
the graph of any Poisson map is coisotropic, and for any Lie subalgebra $\h$ of a Lie algebra $\g$
the annihilator $\h^{\circ}$ is a coisotropic submanifold of the Poisson manifold $\g^*$. Further coisotropic submanifolds $C$ are interesting for a variety of reasons, one being that 
 the distribution $\sharp N^*C$ is a (usually singular) integrable distribution
whose leaf space, if smooth, is a Poisson manifold.

To give a Poisson-analogue of Marle's   result we consider 
\emph{pre-Poisson} submanifolds, i.e. submanifolds $C$ for which $TC+\sharp
N^*C$ has constant rank (or 
equivalently   $ pr_{NC}\circ\sharp \colon N^*C\rightarrow TP|_C \rightarrow NC:=TP|_C/TC$ has constant rank). Natural classes of pre-Poisson submanifolds are given by 
affine subspaces
$\h^{\circ}+\lambda$ of $\g^*$, where $\h$ is a Lie subalgebra of the 
Lie algebra $\g^*$ and $\lambda$ any element of $\g^*$, and of course by coisotropic 
submanifolds and by points.  More details are given  in \cite{CZbis}, where it is also 
shown that  
pre-Poisson submanifolds satisfy some functorial properties. 
This can be used to show that on a  Poisson-Lie group $G$   the graph of $L_h$
 (the left 
translation by some fixed $h\in G$, which clearly is not a Poisson map) is a pre-Poisson
 submanifold, giving rise to a 
  natural constant rank
 distribution $D_h$ on $G$ that leads to interesting constructions.
For instance, if the Poisson structure on $G$ comes from an r-matrix and the point $h$ is
chosen appropriately,
 $G/D_h$ (when  smooth) inherits a
 Poisson structure from $G$, 
and  $[L_h]:G\rightarrow G/D_h$   is a Poisson map
which is moreover equivariant w.r.t. the natural Poisson actions of $G$.

In the following table we characterize submanifolds
of a symplectic or Poisson manifold   in terms of the bundle map $\rho:=pr_{NC}\circ\sharp \colon N^*C\rightarrow  NC$:

\begin{center}
\begin{tabular}{l|l|l}
& $P$ symplectic &$P$ Poisson\\
\hline
$Im(\rho)=0$ & $C$ coisotropic &  $C$ coisotropic\\
$Im(\rho)=NC$ & $C$ symplectic &  $C$ cosymplectic\\
$Rk(\rho)=$const & $C$ presymplectic &  $C$ pre-Poisson\\
\end{tabular}
\end{center}
\vspace*{0.6cm}

In the first part of this paper (sections \ref{embPD}-
 \ref{red}) we consider the Poisson-analog of Marle's result, i.e. we ask
the following question:\begin{quote}
 Given an arbitrary submanifold $C$ of a
Poisson manifold $(P,\Pi)$, under what conditions does there exist
some submanifold $\tilde{P}\subset P$ such that
\begin{itemize}
\item[a)] $\tilde{P}$ has a  Poisson structure induced from $\Pi$
\item[b)] $C$ is a coisotropic submanifold of $\tilde{P}$?
\end{itemize}
When the submanifold $\tilde{P}$ exists, is it unique up to
neighborhood equivalence (i.e. up to a Poisson diffeomorphism on
a tubular neighborhood which fixes $C$)?
\end{quote}
We show in section \ref{embPD} that for any
pre-Poisson submanifold $C$ of a Poisson manifold $P$ there is a
submanifold $\tilde{P}$ which is cosymplectic (and hence  has a
canonically induced Poisson structure) such that $C$ lies
coisotropically in $\tilde{P}$. Further (section \ref{uniqPD})
this cosymplectic submanifold is unique up to neighborhood
equivalence; to the best of our knowledge, this uniqueness result 
is new even in the symplectic setting.
 In section \ref{ex} we give sufficient conditions and necessary
conditions for the existence of a submanifold $\tilde{P}$ as in the above question and we provide examples.  Then in section \ref{red} we deduce statements about
the algebra $C_{bas}^{\infty}(C)$ of functions on $C$ which are basic (invariant),
   meaning that their differentials annihilate
the distribution $\sharp N^*C\cap TC$, and about
and its deformation
quantization. We show that if $C$ is a pre-Poisson submanifold so that the
first and second Lie algebroid cohomology of $N^*C\cap
\sharp^{-1}TC$ vanish, then
 the Poisson algebra of basic functions on $C$ admits a
 deformation quantization.
Finally in section \ref{groids}, assuming that the symplectic
groupoid $\Gamma_s(P)$ of $P$ exists, we describe two subgroupoids
(an isotropic and a presymplectic one) naturally associated to a
pre-Poisson submanifold $C$ of $P$.\\

The second part of this paper (sections \ref{emdpois} and
\ref{aissa}) deals with a different embedding problem, where we
start with an abstract manifold instead of a submanifold of some
Poisson manifold. This is the Poisson-analogue of
Gotay's result.
The question we ask is:
\begin{quote}
 Let $(M,L)$ be a  Dirac manifold. Is there an embedding
$i\colon (M,L)\rightarrow (P,{\Pi})$ into a Poisson manifold such
that
\begin{itemize}
\item[a)] $i(M)$ is a coisotropic submanifold of $P$
\item[b)] the Dirac structure $L$ is induced by the Poisson structure
$\Pi$?
\end{itemize}
Is such an embedding unique up to neighborhood equivalence?
\end{quote}
In the symplectic setting  both existence and uniqueness hold   \cite{Go}.
 One motivation for this question is the deformation
quantization of the Poisson algebra of so-called admissible
functions on $(M,L)$, for a coisotropic embedding as above allows one
to reduce the problem to \cite{CaFeCo2}, i.e. to the deformation
quantization of the basic functions on a coisotropic submanifold
of a Poisson manifold.

It turns out (section \ref{emdpois}) that the above question
admits a positive answer if{f} the distribution $L\cap TM$ on the
Dirac manifold $M$ is regular. In that case one expects the
Poisson manifold $\tP$ to be unique (up to a Poisson
diffeomorphism fixing $M$), provided $\tP$ has minimal dimension.
We are not able to prove this global uniqueness; we can just show
in section \ref{aissa} that the Poisson vector bundle $T\tP|_M$ is
unique (an infinitesimal statement along $M$) and that around each
point of $M$ a small neighborhood of $\tP$ is unique (a local
statement). We remark that A. Wade \cite{Wa} has been considering
a similar question too. Our result about deformation quantization
is the following (Thm. \ref{dqdirac}): let $(M,L)$ be a Dirac
manifold such that $L\cap TM$ has constant rank, and denote by
$\cF$ the regular foliation integrating $L\cap TM$. If the first
and second foliated de~Rham cohomologies of the foliation $\cF$
vanish then the Poisson algebra of admissible functions on $(M,L)$
has a deformation quantization. In Prop. \ref{li} we also notice that the
foliated de~Rham cohomology $\Omega^{\bullet}_{\cF}(M)$
admits the structure of an $L_{\infty}$-algebra (canonically up to 
$L_{\infty}$-isomorphism), generalizing a result of Oh and Park in the presymplectic setting (Thm. 9.4 of \cite{OP}).
\\

We end this introduction describing one of our
   motivations for the first question above, namely an
   application of the Poisson sigma model to quantization
problems.
  The Poisson sigma model is a topological
field theory, whose fields are bundle maps from $T\Sigma$ (for
$\Sigma$ a surface) to the cotangent bundle $T^*P$ of a Poisson
manifold $(P,\Pi)$. It was used by Felder and the first author
\cite{CaFeStar} to derive and interpret Kontsevich's formality
theorem and his star product on the Poisson manifold $P$. The Poisson sigma model with boundary conditions
on a coisotropic submanifold $C$, when suitable assumptions on $C$
are satisfied and $P$ is assumed to be an open subset of $\RR^n$,
provides \cite{CaFeCo1} a deformation quantization of the Poisson
algebra of basic (invariant) functions $C_{bas}^{\infty}(C)$ on
$C$. This result was globalized using a
supergeometric version of Kontsevich's formality theorem
\cite{CaFeCo2}: when the first and second cohomology of the Lie
algebroid $N^*C$ vanish, for $C$ a coisotropic submanifold of any
Poisson manifold $P$, the Poisson algebra $C_{bas}^{\infty}(C)$
admits a deformation quantization. Notice that
  the quotient of $C$ by the distribution $\sharp
N^*C$ is usually not a smooth manifold. Hence
$C_{bas}^{\infty}(C)$ is usually not the algebra of functions on
any Poisson manifold, and
 one cannot apply 
Kontsevich's theorem \cite{K} on deformation quantization of
Poisson manifolds directly.

Calvo and Falceto observed that the most general boundary conditions
for
 the Poisson sigma model  are given by pre-Poisson submanifolds of $(P,\Pi)$ (which they referred to as  ``strongly
regular submanifolds''). They show \cite{CaFa2} that when
$P$ is an open subset of $\RR^n$ the problem of deformation
quantizing the Poisson algebra of basic functions on $C$ can be
reduced to the results of \cite{CaFeCo1}. The computations
in \cite{CaFa2} are carried out choosing local coordinates on $P$
adapted to $C$. The strong regularity condition allows one to choose
local constraints for $C$ such that the number of first class
constraints ($X^{\mu}$s whose Poisson bracket with all other
constrains vanish on $C$) and second class constraints (the
remaining constraints $X^A$, which automatically satisfy
$det\{X^A,X^B\}\neq 0$ on $C$) be constant along $C$. Setting the
second class constraints $X^A$ to zero locally gives a submanifold
with an induced Poisson structure, and the fact that only first
class constraints are left means that $C$ lies in it as a
coisotropic submanifold. Our first question above can be seen as a globalization of Calvo and Falceto's results.\\

\textbf{Conventions:} We use the term
``presymplectic manifold'' to denote a manifold endowed with a
closed 2-form of \emph{constant rank}, i.e. such that its kernel have
constant rank. However we stick to the denominations ``presymplectic groupoid'' coined in
\cite{BCWZ} and ``presymplectic leaves'' (of a Dirac manifold) despite the fact that the 2-forms
on these objects do not have constant rank, for these denominations seem to be established
in the literature.
\\

\textbf{Acknowledgements:} As M.Z. was a graduate student
Marius Crainic first called to his attention some of the questions
discussed in section \ref{emdpois}, and
some of the existence results obtained
in  the same section arose from discussion between
M.Z. and
 Alan Weinstein, who at the time was his thesis advisor and whom
 he gratefully thanks.
  M.Z. also would
like to thank Aissa Wade for pointing out the necessity of a
minimal dimension condition mentioned in Section \ref{aissa} and
Eva Miranda for showing him the reference \cite{AM}. We also thank Rui Loja Fernandes 
for comments and the referees    
for suggesting valuable improvements to a 
previous version of this manuscript.
A.S.C. acknowledges partial support of SNF Grant No.~20-113439.
This work has been partially supported
by the European Union through the FP6 Marie Curie RTN ENIGMA (Contract
number MRTN-CT-2004-5652) and by the European Science Foundation
through the MISGAM program. M.Z.   acknowledges support from the Forschungskredit 2006 of the University of Z\"urich.

\section{Basic definitions}\label{def}

We will use some notions from Dirac linear algebra \cite{Cou} \cite{BR}. A Dirac
structure on a vector space $P$ is a subspace $L\subset P\oplus
P^*$ which is maximal isotropic w.r.t. the natural symmetric inner
product on $P\oplus P^*$ (i.e. $L$ is isotropic and has same
dimension as $P$). A Dirac structure $L$ specifies a subspace
$\cO$, defined as the image of $L$ under the projection $P\oplus
P^*\rightarrow P$, and a skew-symmetric bilinear form $\omega$ on
$\cO$, given by $\omega(X_1,X_2)=\langle \xi_1,X_2 \rangle$ where
$\xi_1$ is any element of $P^*$ such that $(X_1,\xi_1)\in L$. The
kernel of $\omega$ (which in terms of $L$ is given as $L\cap P$)
is called \emph{characteristic subspace}. Conversely, any choice
of bilinear form defined on a subspace of $P$ determines a Dirac
structure on $P$. Given this equivalence, we will
sometimes work with the bilinear form $\omega$ on $\cO$ instead of
working with $L$.

We consider now Poisson vector spaces $(P,\Pi)$ (i.e. $\Pi\in
\wedge^2P$; we denote by $\sharp\colon P^*\rightarrow P$  the map
induced by contraction with $\Pi$). The Poisson structure on $P$
is encoded by the Dirac structure $L_P=\{(\sharp \xi,\xi)\colon
\xi\in P^*\}$. The image of $L_P$ under the projection onto the
first factor is $\cO=\sharp P^*$, and the bilinear form $\omega$
is non-degenerate.

\begin{remark}\label{intr}
We recall that any subspace $W$ of a Dirac vector space $(P,L)$
has an induced Dirac structure $L_W$; the bilinear form
characterizing $L_W$ is just the pullback of $\omega$ (hence it is
defined on $W\cap \cO$). When $(P,\Pi)$ is actually a Poisson vector
space, one shows 
 that the symplectic orthogonal of $W\cap \cO$ in
$(\cO,\omega)$ is $\sharp W^{\circ}$. Hence $\sharp W^{\circ}\cap W$
 is the kernel of the
restriction of $\omega$ to $W\cap \cO$, i.e. it is the characteristic subspace of
 the Dirac structure $L_W$, and we will refer to it as the
 \emph{characteristic subspace of W}.
Notice that pulling back Dirac structure is functorial \cite{BR}
(i.e. if $W$ is contained in some other subspace $W'$ of $P$,
pulling back $L$ first to $W'$ and then to $W$ gives the Dirac
structure $L_W$), hence $L_W$, along with the corresponding
bilinear form and characteristic subspace, is \emph{intrinsic} to $W$.\\
\end{remark}

Let $W$ be a subspace of the Poisson vector space $(P,\Pi)$.
 $W$ is
called \emph{coisotropic} if $ \sharp W^{\circ}\subset W$, which
by the above means that $W\cap \cO$ is coisotropic in
$(\cO,\omega)$.

 $W$ is called
\emph{Poisson-Dirac} subspace \cite{CrF} when  $\sharp
W^{\circ}\cap W=\{0\}$; equivalent conditions are that
  $W\cap \cO$ be a symplectic subspace of
$(\cO,\omega)$ or that
  the pullback Dirac structure $L_{P}$ correspond to a
Poisson bivector on $W$.  The Poisson bivector on $W$
is described as follows \cite{CrF}: its sharp map $\sharp_W \colon W^*\rightarrow W$
is given by $\sharp_W \tilde{\xi}=\sharp \xi$, where $\xi\in P^*$
 is any extension of $\tilde{\xi}$ which annihilates $\sharp W^{\circ}$.

$W$ is called \emph{cosymplectic}
subspace if $\sharp W^{\circ}\oplus W=P$, or equivalently if the
pushforward of $\Pi$ via the projection $P\rightarrow P/W$ is an
invertible bivector. Notice that if $W$ is cosymplectic then it
has a canonical complement $\sharp W^{\circ}$ which is  a
symplectic subspace of $(\cO,\omega)$. Clearly a cosymplectic
subspace is automatically a Poisson-Dirac subspace.\\

Now we pass to the global definitions.
A Dirac structure on $P$ is a maximal isotropic subbundle $L\subset TP\oplus T^*P$
which is integrable, in the sense that its sections are closed under the so-called Courant bracket (see \cite{Cou}). The image of $L$ under
the projection onto the first factor is an integrable singular distribution, whose leaves
(which are called presymplectic leaves) are endowed with closed 2-forms. A Poisson structure
on $P$ is a bivector $\Pi$ such that $[\Pi,\Pi]=0$.

Coisotropic
and cosymplectic submanifolds of a Poisson manifold are defined
exactly as in the linear case; a Poisson-Dirac submanifold
additionally requires that the bivector induced on the submanifold
by the point-wise condition be smooth \cite{CrF}. Cosymplectic
submanifolds  
 are automatically Poisson-Dirac
submanifolds (the smoothness of the induced bivector is ensured because
 $L_P\cap (\{0\}\oplus N^*\tP)$ has constant rank zero).
The Poisson bracket on a Poisson-Dirac  submanifold $\tP$ of $(P,\Pi)$ is computed as follows:
$\{\tilde{f}_2,\tilde{f}_2\}_{\tP}$ is the restriction to $\tP$ of $\{f_1,f_2\}$, where the $f_i$ are extensions of $\tilde{f}_i$ to $P$
such that $df_i|_{\sharp N^*\tP}=0$ (for at least one of the two functions).
We will also need a definition which does not have a linear algebra counterpart:
\begin{defi}
A submanifold $C$ of a Poisson manifold $(P,\Pi)$ is called
\emph{pre-Poisson} if the rank of $TC+\sharp N^*C$ is constant
along $C$.
\end{defi}

\begin{remark}
An alternative characterization of pre-Poisson submanifolds is the requirement that
  $\Pi|_{\wedge^2 N^*C}$  (or equivalently the corresponding sharp map $ pr_{NC}\circ\sharp \colon N^*C\rightarrow TP|_C \rightarrow NC:=TP|_C/TC$) have constant rank.
 Indeed the kernel of  
$N^*C\rightarrow NC$ is $N^*C\cap \sharp^{-1}TC$, which is the annihilator of 
$TC+\sharp N^*C$. The map $N^*C\rightarrow NC$ is identically zero iff $C$ is coisotropic and is an isomorphism iff $C$ is cosymplectic.
\end{remark}

Calvo and Falceto already considered \cite{CaFa1}\cite{CaFa2} such
submanifolds and called them ``strongly regular submanifolds''.
We prefer to call them ``pre-Poisson'' because when $P$ is a symplectic manifold they
reduce to presymplectic submanifolds\footnote{Further reasons are the following:
the subgroupoid associated to a pre-Poisson
manifold, when it exists, is presymplectic (see Prop. \ref{pre}). The Hamiltonian version of the
Poisson Sigma Model with boundary conditions on $P$ (at $t=0$) and on a submanifold $C$ (at $t=1$)
delivers a space of solutions which is presymplectic if{f} $C$ is pre-Poisson.}.
See Section \ref{ex} for
several examples.

\section{Existence of coisotropic embeddings for pre-Poisson submanifolds}\label{embPD}

In this section we consider the problem of embedding a submanifold of a Poisson manifold
coisotropically in a Poisson-Dirac submanifold, and show that this can be always done
for pre-Poisson submanifolds.

We start with some linear algebra.
\begin{lemma}\label{poisLA}
Let $(P,\Pi)$ be a Poisson vector space and $C$ a subspace. The
Poisson-Dirac subspaces of $P$ in which $C$ sits coisotropically
are exactly the subspaces $W$
satisfying
\begin{eqnarray}\label{condPoisson}
W+\sharp C^{\circ}\supset \cO\\
\label{condint} W \cap (C+\sharp C^{\circ})=C
\end{eqnarray}
where $\cO=\sharp P^*$.
Among the  Poisson-Dirac subspaces above the cosymplectic ones are
exactly those of maximal dimension, i.e. those for which
$W+\sharp C^{\circ}=P$.
\end{lemma}

\begin{remark}\label{convenient}
It is often more convenient to work with the following characterization of the Poisson-Dirac subspaces $W$ containing coisotropically $C$: $W=R\oplus C$, where the subspace $R$ satisfies
\begin{eqnarray}\label{condPois}
R\oplus(C+\sharp C^{\circ})\supset \cO.
\end{eqnarray}
Among these, the cosymplectic subspace are those for which $R$ satisfies the stronger condition $R\oplus(C+\sharp C^{\circ})=P$. 
When $\Pi$ corresponds to a linear symplectic form $\omega$, both conditions become
$R\oplus(C+C^{\omega})=P$.

\end{remark}

\begin{proof}
The condition that $W$ be a Poisson-Dirac subspace is 
\begin{equation}\label{cond1}
W \cap \sharp W^{\circ}=0.
\end{equation} Let us denote by $\sharp_W$ the sharp map of the induced bivector on $W$.
The condition that $C$ is contained in $W$ coisotropically 
is $\sharp_W \tilde{\xi} \in C$ for all elements $\tilde{\xi} \in W^*$ annihilating $C$. $\sharp_W \tilde{\xi}$ 
is obtained extending $\tilde{\xi}$ to some $\xi\in (\sharp W^{\circ})^{\circ}=\sharp^{-1}W$ and applying   $\sharp$. Hence the condition that $C$ is contained in $W$ coisotropically 
can be phrased as
\begin{equation}\label{cond2}
W \cap \sharp C^{\circ}\subset C\subset W.
\end{equation} 
 
We show now that conditions \eqref{cond1} and \eqref{cond2} are equivalent to conditions \eqref{condPoisson} and \eqref{condint}.

 We  have $\eqref{cond2} \Rightarrow  \eqref{condint}$, because due to $C\subset W$ 
 we have $W \cap (C+\sharp C^{\circ})=C+(W\cap \sharp C^{\circ})$. The implication
 $\eqref{cond2} \Leftarrow  \eqref{condint}$ is immediate.
 
Now assume that either of  \eqref{cond2} or  \eqref{condint} hold true. 
Applying $\sharp(\bullet)^{\circ}$ we see that  condition \eqref{condPoisson} is equivalent to $C \cap \sharp W^{\circ}=\{0\}$.
Since applying condition \eqref{cond2} we have $$W \cap \sharp W^{\circ}=
(W \cap \sharp C^{\circ})\cap \sharp W^{\circ}\subset
 C \cap \sharp W^{\circ}\subset W \cap \sharp W^{\circ},$$
  the equivalence of conditions \eqref{cond1} and \eqref{condPoisson} is proven.

To prove the last statement of the lemma let $W$ satisfy eq. \eqref{condPoisson} and \eqref{condint}; in particular
$W$ is Poisson-Dirac. By dimension counting $W$ is
cosymplectic if{f} the restriction of $\sharp$ to $W^{\circ}$ is
injective, i.e. if{f} $W^{\circ}\cap \cO^{\circ}=\{0\}$ or
$W+\cO=P$. Using eq. \eqref{condPoisson}this is seen to be equivalent to $W+\sharp C^{\circ}=P$.\end{proof}

Now we pass from linear algebra to global geometry. Given a submanifold $C$ of a Poisson manifold $P$, one might try to construct a Poisson-Dirac submanifold in which $C$ embeds coisotropically
applying the corresponding symplectic construction ``leaf by leaf'' in a smooth way. It would then be natural to require that the characteristic ``distribution'' $TC\cap \sharp N^*C$ 
of $C$ have constant rank. However this approach generally does not work, because even when it has constant rank
 $TC\cap \sharp N^*C$ might not be smooth (see example \ref{noncont}). 
The right condition to ask is instead that $TC+\sharp N^*C$ have constant rank:

\begin{thm}\label{emb}
Let $C$ be a pre-Poisson submanifold of a Poisson manifold
$(P,\Pi)$. Then there exists a
cosymplectic submanifold $\tP$ containing $C$ such that $C$ is
coisotropic in $\tP$.
\end{thm}
\begin{proof}
Because of the rank condition on $C$ we can choose a smooth
subbundle $R$ of $TP|_C$ which is a complement  to $TC + \sharp
N^*C$.  Then by Lemma \ref{poisLA} at every point $p$ of $C$ we
have that $T_pC\oplus R_p$ is a cosymplectic subspace of $T_pP$ in
which $T_pC$ sits coisotropically. ``Thicken''  $C$ to a smooth
submanifold  $\tP$ of $P$ satisfying $T\tP|_C=TC\oplus R$.  
Since 
$T_p\tP\oplus
\sharp N_p^*\tP=T_pP$ is an open condition that holds at every point $p$ of $C$, it holds at points in a tubular neighborhood of $C$ in $\tP$. Hence, shrinking $\tP$ if necessary, we obtain a cosymplectic submanifold of $P$ containing  coisotropically $C$.
\end{proof}

\begin{remark}\label{nons}
 The above proposition says that if $C$ is a pre-Poisson submanifold then we can
choose a subbundle $R$ over $C$ with fibers as in eq. \eqref{condPois}
and
 ``extend`` $C$ in direction of $R$ to obtain a
Poisson-Dirac  submanifold of $P$ containing $C$ coisotropically.
If $C$ is not a pre-Poisson submanifold of $(P,\Pi)$, we might
still be able to find a smooth bundle $R$ over $C$ consisting of
subspaces
 as in eq. \eqref{condPois}. However ``extending'' $C$ in direction of this
subbundle will usually not give a submanifold with a smooth
Poisson-Dirac structure, see Example \ref{mariusrui} below.
\end{remark}

Now we deduce consequences about Lie algebroids. See section
\ref{groids} for the corresponding ``integrated'' statements.

\begin{lemma}\label{coro}
Let $C$ be a subspace of a Poisson vector space $(P,\Pi)$ and $W$ a cosymplectic subspace containing $C$ as a coisotropic subspace.  Then
$C+\sh C^{\circ}=C\oplus \sh W^{\circ}$.
\end{lemma}
\begin{proof}   
 The
inclusion ``$\supset$'' holds because $C\subset W$. The
other inclusion follows   by this
argument: write any $\xi\in C^{\circ}$ uniquely as $\xi_1+\xi_2$ where
$\xi_1$ annihilates $\sh W^{\circ}$ and $\xi_2$ annihilates
$W$. Then $\sharp \xi_1=\sharp_W(\xi_1|_{W})\in C$,
where ${\sharp_W}$ denotes the sharp map of $W$, since $C$
is coisotropic in $W$. Hence $\sharp \xi=\sharp
\xi_1+\sharp\xi_2\in C+ \sh W^{\circ}$. Finally, we have a direct sum in $ C\oplus \sh W^{\circ}$
because $\sh W^{\circ} \cap W=\{0\}$ and $C\subset W$.
\end{proof}

\begin{prop}\label{algoid}
Let $C$  be a submanifold of a Poisson manifold $(P,\Pi)$. Then
$N^*C\cap \sharp^{-1}TC$ is a Lie subalgebroid
  of $T^*P$ if{f} $C$ is pre-Poisson. In that case, for
any cosymplectic submanifold $\tP$ in which $C$ sits
coisotropically, $N^*C\cap \sharp^{-1}TC$ is isomorphic as a Lie algebroid to the
annihilator of $C$ in $\tP$.
\end{prop}
\begin{proof}
At every point $N^*C\cap \sharp^{-1}TC$ is the annihilator of
$TC+\sharp N^*C$, so it is a vector bundle if{f} $C$ is
pre-Poisson. So assume that $C$ be pre-Poisson.
For any
cosymplectic submanifold $\tP$ 
the   embedding $T^*\tP \rightarrow T^*P$,
obtained  extending a covector in $T^*\tP$ so that it
annihilates $\sharp N^*\tP$, is a Lie algebroid morphism (Cor.
2.11 and Thm. 2.3 of \cite{Xu}). If  $C$ lies coisotropically in $\tilde{P}$, 
by Lemma \ref{coro} $TC+\sh
N^*C=TC\oplus \sh N^*\tP|_C$. Hence $N^*_{\tP}C$, the conormal bundle
of $C$ in $\tP$, is mapped \emph{isomorphically} onto $(TC\oplus
\sh N^*\tP)^{\circ}=(TC + \sharp N^*C)^{\circ}=N^*C\cap
\sharp^{-1}TC$. Since $N^*_{\tP}C$ is a Lie subalgebroid of
$T^*\tP$ \cite{Ca}, we are done.
\end{proof}

\begin{remark}
The fact that $N^*C\cap \sharp^{-1}TC$ is a Lie algebroid if $C$
is pre-Poisson can also be deduced as follows. The Lie algebra
$(\cF\cap I)/I^2$ forms a Lie-Rinehart algebra over the
commutative algebra $C^{\infty}(P)/I$, where  $I$ is the vanishing
ideal of $C$ and $\cF$ its Poisson-normalizer in $C^{\infty}(P)$.
Lemma 1 of \cite{CaFa1} states that $C$ being pre-Poisson is
equivalent to
  $N^*C\cap \sharp^{-1}TC$ being spanned by
differentials of functions in $\cF\cap I$. From this one deduces
easily that $(\cF\cap I)/I^2$ is  identified with the sections of
  $N^*C\cap \sharp^{-1}TC$, and since $C^{\infty}(P)/I$ are just
  the smooth functions on $C$ we deduce that   $N^*C\cap
  \sharp^{-1}TC$ is a Lie algebroid over $C$.
\end{remark}

\section{Uniqueness of coisotropic embeddings for pre-Poisson
submanifolds}\label{uniqPD}

Given a submanifold $C$ of a Poisson manifold $(P,\Pi)$ in this
section we investigate the uniqueness (up Poisson diffeomorphisms fixing $C$)
of \emph{cosymplectic}
submanifolds in which $C$ is embedded coisotropically.

This lemma tells us that we need to consider only the case in which  $C$
is pre-Poisson and the construction of Thm. \ref{emb}:
\begin{lemma}
A submanifold $C$ of a Poisson manifold $(P,\Pi)$ can be embedded
coisotropically in a cosymplectic submanifold $\tP$ if{f} it is
pre-Poisson. In this case all such $\tP$ are constructed (in a neighborhood of $C$)
as in Thm. \ref{emb}.
\end{lemma}
\begin{proof}
In Thm. \ref{emb} we saw that given any pre-Poisson submanifold
 $C$,   choosing a smooth
 subbundle $R$ with $R\oplus(TC + \sharp N^*C)=TP|_{C}$ and
 ``thickening'' $C$ in direction of $R$ gives a submanifold $\tP$
 with the required properties.

Now let $C$ be any submanifold embedded coisotropically in
 a cosymplectic submanifold $\tP$. By Remark \ref{convenient}, for any complement $R$ of $TC$ in
 $T\tP|_{C}$ we have $R\oplus(TC + \sharp N^*C)=TP|_{C}$. This has two consequences:
 first
 the rank of $TC + \sharp N^*C$ must be constant, concluding the
 proof of the ``if{f}'' statement of the  lemma. Second, it proves
 the final statement of the lemma.
\end{proof}

  When $C$ is a point $\{p\}$ then
$\tP$ as above is a slice transverse to the symplectic leaf through $p$
(see Ex. \ref{exstrong}) and $\tP$ is  unique up Poisson diffeomorphism by Weinstein's splitting theorem (Lemma 2.2 in
\cite{We}; see also  Thm. 2.16 in \cite{Va}). A generalization of
its proof gives

\begin{prop}\label{uniq}
Let $\tP_0$ be a cosymplectic submanifold of a Poisson manifold
$P$ and $\pi\colon U\rightarrow\tP_0$ a projection of some
tubular neighborhood of $\tP_0$ onto $\tP_0$. Let $\tP_t$,
$t\in\RR$, be a smooth family of cosymplectic submanifolds such
that all $\tP_t$ are images of sections of $\pi$. Then, for $t$
close enough to zero, there are Poisson diffeomorphisms $\phi_t$ mapping
open sets of $\tP_0$ to open sets of  $\tP_t$. The $\phi_t$'s can be chosen so that 
 the curves
$t\mapsto \phi_t(y)$ (for $y\in \tP_0$) are tangent to $\sharp
N^*\tP_{t}$ at time $t$. 
\end{prop}
\begin{proof}
We will use the following fact, whose straightforward proof we
omit: let $\tP_t$, $t\in \RR$, be a smooth family of submanifolds
of a manifold $U$, and $Y_t$ a time-dependent vector field on $U$.
Then $Y+\frac{\partial}{\partial t}$ (considered as a vector field on $U\times \RR$) is
tangent to the submanifold $\bigcup_{t\in\RR}(\tP_t,t)$  if{f}
for each $\bar{t}$ and each integral curve $\gamma$ of $Y_t$ in $U$ with
$\gamma(\bar{t})\in \tP_{\bar{t}}$ we have $\gamma(t)\in \tP_t$
(at all times where $\gamma$ is  defined).

Denote by $s_t$ the section of $\pi$ whose image is $\tP_t$. We
are interested in time-dependent vector fields $Y_t$ on $U$
such that for all $\bar{t}$ and $y\in \tP_{\bar{t}}$
\begin{equation}\label{vfcond}
Y_{\bar{t}}(y)={s_{\bar{t}}}_*(\pi_*Y_y)+\frac{d}{dt}|_{\bar{t}}s_t(\pi(y)).
\end{equation}
We claim that, for such a vector field,
$(Y+\frac{\partial}{\partial t})$ is tangent to
$\bigcup_{t\in\RR}(\tP_t,t)$. Indeed
\begin{eqnarray}
(Y+\frac{\partial}{\partial t})(y,\bar{t})&=&Y_{\bar{t}}(y)+\frac{\partial}{\partial t}\\
&=&{s_{\bar{t}}}_*(\pi_*
Y_y)+\frac{d}{dt}|_{\bar{t}}s_t(\pi(y))+\frac{\partial}{\partial
t}.
\end{eqnarray}
Since ${s_{\bar{t}}}_*(\pi_* Y_y)$ is tangent to
$(\tP_{\bar{t}},\bar{t})$, and
$\frac{d}{dt}|_{\bar{t}}s_t(\pi(y))+\frac{\partial}{\partial t}$
is the velocity at time $\bar{t}$ of the curve $(s_t(\pi(y)),t)$,
the claimed tangency follows. Hence by the fact recalled in the
first paragraph we deduce that the flow $\phi_t$ of $Y_t$ takes
points $y$ of $\tP_0$ to  $\tP_{\bar{t}}$ (if $\phi_t(y)$ is defined until time
$\bar{t})$.

So we are done if we realize such  $Y_t$ as the hamiltonian vector
fields of a smooth family of functions $H_t$ on $U$. For each
fixed $\bar{t}$, eq. \eqref{vfcond} for $Y_{\bar{t}}$ is just a
condition on the second  component of
$Y_{\bar{t}}\in T_yP=T_y\tP_{\bar{t}}\oplus ker_y\pi_*$ for all $y\in \tP_{\bar{t}}$, and the second component is
determined exactly by the action of $Y_{\bar{t}}$ on functions $f$
vanishing on $\tP_{\bar{t}}$. We have
$$Y_{\bar{t}}(f)=X_{H_{\bar{t}}}(f)=-dH_{\bar{t}}(\sharp df),$$ and the restriction of $\sharp$ to
$N^*\tP_{\bar{t}}$ is injective because  $\tP_{\bar{t}}$
is cosymplectic. Together we obtain that
specifying the vertical component of $X_{H_{\bar{t}}}$ at points
of $\tP_{\bar{t}}$ is equivalent to specifying the derivative of
$H_{\bar{t}}$ in direction of $\sharp N^*\tP_{\bar{t}}$, which is
transverse to $\tP_{\bar{t}}$. We can clearly find a function
$H_{\bar{t}}$ satisfying the required conditions on its derivative
along $\tP_{\bar{t}}$, i.e. so that $X_{H_{\bar{t}}}$ satisfies \eqref{vfcond}.
 Choosing $H_t$ smoothly for every $t$ we
conclude that the flow $\phi_t$ of $X_{H_{t}}$, which obviously consists
of Poisson diffeomorphisms, will take $\tP_0$ (or rather any
subset of it on which the flow is defined up to time $\bar{t}$) to
$\tP_{\bar{t}}$.

Choosing each $H_t$ so that it vanishes on
$\tP_t$ delivers a flow $\phi_t$ ``tangent'' to the $\sharp N^*\tP_{t}$'s.
\end{proof}

Now we are ready to prove the  uniqueness of $\tP$:
\begin{thm}\label{uniPhi}
Let $C$ be a pre-Poisson submanifold $(P,\Pi)$, and $\tP_0$,
$\tP_1$ cosymplectic submanifolds that contain $C$ as a coisotropic submanifold. Then,
shrinking $\tP_0$ and $\tP_1$ to a smaller tubular
neighborhood of $C$ if necessary, there is a Poisson diffeomorphism $\Phi$ from
$\tP_0$ to $\tP_1$ which is the identity on $C$.
\end{thm}

\begin{proof}
In a neighborhood $U$ of $\tP_0$ take a projection $\pi\colon
U\rightarrow\tP_0$; choose it so that at points of $C\subset \tP_0$ the fibers
of $\pi$ are tangent to $\sharp N^*\tP_0|_C$. For $i=0,1$ make some choices of
maximal dimensional subbundles $R_i$ satisfying eq. \eqref{condPois}
to write $T\tP_i|_C=TC\oplus R_i$, and choose a smooth curve
of subbundles $R_t$ satisfying eq. \eqref{condPois} and agreeing with
$R_0$ and $R_1$ at $t=0,1$
(there is no
topological obstruction to this because $R_0$ and $R_1$ are both
complements to the same subbundle $TC+\sharp N^*C$). By Thm.
\ref{emb} we obtain a curve of cosymplectic submanifolds $\tP_t$,
which  moreover by Lemma \ref{coro}
at points of $C$ are all transverse to $\sharp N^*\tP_0|_C$, i.e. to
the fibers of $\pi$.

Hence we  are in the situation of Prop. \ref{uniq}, which allows us
to construct a Poisson diffeomorphism from $\tP_0$ to $\tP_t$ for
small $t$.  Since $C\subset \tP_t$ for all $t$, in the proof of
Prop.\ref{uniq} we have that the sections $s_t$ are trivial on
$C$, hence by eq. \eqref{vfcond} the second component of $X_{H_t}\in
T_y\tP_{{t}}\oplus ker_y\pi_*$ at points $y$ of
$C\subset \tP_t$ is zero. Choosing $H_t$ to vanish on $\tP_t$ we
obtain $X_{H_t}=0$ at points of $C\subset \tP_t$. From this we
deduce two things: in a tubular neighborhood of $C$ the flow
$\phi_t$ of $X_{H_t}$ is defined for all $t\in [0,1]$, and each
$\phi_t$ keeps points of $C$ fixed. Now just let $\Phi:=\phi_1$.
\end{proof}

The derivative at points of $C$ of the Poisson diffeomorphism
$\Phi$ constructed in Thm. \ref{uniPhi} gives an isomorphism of Poisson vector
bundles $T{\tP_0}|_C\rightarrow T{\tP_1}|_C$ which is the identity
on $TC$. The construction of $\Phi$ involves many
choices; we wish now to give a \emph{canonical} construction for such a vector bundle isomorphism. We first need a linear algebra lemma.

\begin{lemma}\label{infuLA}
Let $C$ be a subspace of a Poisson vector space $(P,\Pi)$ and $V,W$ two cosymplectic subspaces containing $C$ as a coisotropic subspace. There exists a \emph{canonical}
isomorphism of Poisson vector spaces
 $\varphi\colon V\rightarrow W$ which is the identity on $C$.
\end{lemma}

\begin{proof}Notice that $V$ and $W$ have the same dimension by Lemma \ref{poisLA}.
 First
we consider  
$$A\colon V\rightarrow \sharp V^{\circ}$$
determined by the requirement that $W=\{v+Av:v\in V\}$. $A$
is well-defined since $ \sharp V^{\circ}$ is a complement in $P$ both
to $V$ (because $V$ is cosymplectic) and to $W$ (because
$W\cap (C+\sharp C^{\circ})=C$ by Lemma \ref{poisLA} and $C+\sh
C^{\circ}=C\oplus \sharp V^{\circ}$ by Lemma \ref{coro}). Notice that, since
$C$ lies in both $V$ and $W$, the restriction of $A$ to $C$
is zero.

Now we mimic a construction
in symplectic linear algebra \cite{Ana} where one deforms canonically a
complement of a coisotropic subspace $C$ to obtain an isotropic
complement. 
 We deform $A+Id$ by adding
$$B\colon V \rightarrow C\cap \sharp C^{\circ}, v\mapsto
\frac{1}{2}{\sharp_V}(\Omega(Av,A\bullet)).$$ Here
$\sharp_V$ is the sharp map of the cosymplectic submanifold
$V$ and $\Omega$ denotes the symplectic form on $\cO:=\sharp P^{*}$.
  $B$ is
well-defined  because the element  $\Omega(Av,A\bullet)$ of $V^*$
annihilates $C$ (recall that $A|_{C}=0$) and because $C$ is
coisotropic in $V$. Further it is clear that the restriction of
$B$ to $C$ is zero.

At this point we are ready to define
$$\varphi\colon V \rightarrow W, v\mapsto v+Av+Bv.$$
This is  well-defined  (since $C\cap \sharp C^{\circ}\subset W$) and 
 is an isomorphism: if $v+Bv+Av=0$ then
$v+Bv=0 $ and $Av=0$ (because $V$ is transversal to $\sharp
V^{\circ}$); from $Av=0$ we deduce $Bv=0$ hence $v=0$.
To show that $\varphi$ matches the linear Poisson structures on $V$ and $W$ we notice that 
   $\varphi$
restricts to a map from $V\cap\cO$ to
$W\cap \cO$ (because the images if $A$ and $B$ lie in
$\cO$). This restriction  is  an isomorphism 
because
source and target have the same dimension 
(they both contain
$C\cap \cO$ as a coisotropic subspace); 
we show that it is a linear
symplectomorphism. If $v_1,v_2\in V\cap\cO$ we have
$\Omega(\varphi v_1,\varphi
v_2)=\Omega(v_1+Bv_1,v_2+Bv_2)+\Omega(Av_1,Av_2)$, for the cross
terms vanish since $A$ takes values in $\sharp V^{\circ}$. Now
$\Omega(Bv_1,\bullet)|_{V\cap \cO}=-\frac{1}{2}\Omega(Av_1,A\bullet)|_{V\cap\cO}$
using the fact that $\Omega(\sharp \xi,\bullet)=-\xi|_{\cO}$  for
any covector $\xi$ of $P$. Further $\Omega(Bv_1,Bv_2)$ vanishes
because $B$ takes values in $C\cap \sharp C^{\circ}$. So
altogether we obtain $\Omega(\varphi v_1,\varphi v_2)=\Omega(v_1,
v_2)$ as desired.
\end{proof}

\begin{prop}\label{infu}
Let $C$ be a pre-Poisson submanifold $(P,\Pi)$, and $\tP$, $\hP$
cosymplectic submanifolds that contain $C$ as a coisotropic submanifold.
Then there is a \emph{canonical}
isomorphism of Poisson vector bundles
 $\varphi\colon T\tP|_C\rightarrow T\hP|_C$ which is the identity on $TC$.
\end{prop}
\begin{proof}
At each point $p\in C$ we construct $\varphi_p$  
applying Lemma \ref{infuLA} to $V=T_p \tP$ and $W=T_p\hP$.
We want to check  that the resulting map $\varphi\colon T\tP|_C\rightarrow T\hP|_C$ is smooth (this is not clear a priori because the construction of Lemma \ref{infuLA} involves the symplectic leaves $\cO$ of $P$, which  may be of different dimensions). 
It is enough to check that  
 if $X$ is a smooth section of $\sharp N^*\tP|_C$, then
$\Omega(X,\bullet)|_{\sharp N^*\tP}:\sharp N^*\tP\rightarrow \RR$
is smooth. This follows from the fact that $\tP$ is
cosymplectic: since $\sharp\colon  N^*\tP\rightarrow \sharp
N^*\tP$ is bijective, there is a smooth section $\xi$ of $N^*\tP$
with $\sharp \xi=X$, and $\Omega(X,\bullet)|_{\sharp
N^*\tP}=\xi|_{\sharp N^*\tP}$. Altogether we obtain that 
$\varphi$ is a smooth,  canonical
isomorphism of Poisson vector bundles.
\end{proof}

\begin{remark} The isomorphism $\varphi\colon T\tP|_C\rightarrow T\hP|_C$ constructed in Prop. \ref{infu}
can be extended to a Poisson vector bundle automorphism  of
$TP|_C$, by applying the following at each point of $C$.

The linear isomorphism $\varphi\colon V \rightarrow W$ of Lemma \ref{infuLA} (using the notation of the lemma)
can be extended to a Poisson  automorphism  of
$P$ as follows: define
$$(\varphi,pr)\colon  V\oplus \sharp V^{\circ} \rightarrow W\oplus \sharp W^{\circ}$$ where $pr$ denotes the projection of $\sharp V^{\circ}$ onto
$\sharp W^{\circ}$ along $C$ (recall from Lemma \ref{coro} that $C\oplus
\sharp V^{\circ}=C\oplus \sharp W^{\circ}$). $(\varphi,pr)$ restricts to a linear
automorphism of $\cO=(V\cap \cO)\oplus \sharp V^{\circ}$
 which preserves the symplectic form: the
only non-trivial check is $\Omega(pr (v_1),pr
(v_2))=\Omega(v_1,v_2)$ for $v_i \in \sharp V^{\circ}$, which follows
because $pr (v_i)-v_i\in C\cap \sharp C^{\circ}$.

\end{remark}

\begin{remark} The isomorphism $\varphi$ constructed in Prop. \ref{infu}
can be extended to a Poisson vector bundle automorphism  of
$TP|_C$ as follows: define
$$(\varphi,pr)\colon  T\tP\oplus \sharp N^*\tP \rightarrow T\hP\oplus \sharp
N^*\hP$$ where $pr$ denotes the projection of $N^*\tP$ onto
$N^*\hP$ along $TC$ (recall from Lemma \ref{coro} that $TC\oplus
N^*\tP=TC\oplus N^*\hP$). $(\varphi,pr)$ restricts to a linear
automorphism of $T\cO=(T\tP\cap T\cO)\oplus \sharp N^*\tP$
 which preserves the symplectic form: the
only non-trivial check is $\Omega(pr (v_1),pr
(v_2))=\Omega(v_1,v_2)$ for $v_i \in \sharp N^*\tP$, which follows
because $pr (v_1)-v_1\in TC\cap \sharp N^*C$.
\end{remark}

\section{Conditions and examples}\label{ex}

Let $C$ be as usual a submanifold of the Poisson manifold
$(P,\Pi)$; in Section \ref{embPD} we considered the question of
existence of a Poisson-Dirac submanifold $\tP$ of $P$ in which $C$
is contained coisotropically. In Thm. \ref{emb} we showed that a
\emph{sufficient} condition is that $C$ be pre-Poisson, which by
Prop. \ref{algoid} is equivalent to saying that $N^*C\cap
\sharp^{-1}TC$ be a Lie algebroid.

A \emph{necessary} condition is that the (intrinsically defined)
characteristic distribution $TC \cap \sharp N^*C$ of $C$ be the
distribution associated to a Lie algebroid over $C$; in particular
its rank locally can only increase. This is a necessary condition
since the concept of characteristic distribution is an intrinsic
one (see Remark \ref{intr}), and the characteristic distribution
of a coisotropic submanifold of a Poisson manifold is the image of
the anchor of its conormal bundle, which is a Lie algebroid.
 
The submanifolds $C$ which are not covered by the above conditions
are those for which $N^*C\cap \sharp^{-1}TC$ is not a Lie
algebroid but its image $TC \cap \sharp N^*C$ under $\sharp$ is
the image of
the anchor of some Lie algebroid over $C$. Diagrammatically:\\

\begin{itemize}
\item[] $\{ C$ s.t. $N^*C\cap \sharp^{-1}TC$
 is a Lie algebroid, i.e. $C$ is pre-Poisson $\} \subset$\\
\item[] $\{C $ sitting coisotropically in some
Poisson-Dirac
submanifold $\tilde{P}$ of $P \} \subset$\\
\item[] $\{ C$ s.t. $TC \cap \sharp N^*C$
 is the distribution of some Lie algebroid over  $C\}.$\\
\end{itemize}

In the remainder of this section we present examples of the above situations.
We start with basic examples of pre-Poisson submanifolds;  we refer the reader to Section 6 of \cite{CZbis} for examples  in which the Poisson manifold $P$ is
the dual of a Lie algebra and $C$ an affine subspace.

\begin{ep}\label{exstrong}
An obvious example is when $C$ is a coisotropic submanifold of
$P$,  and in this
case the construction of Thm. \ref{emb} delivers $\tP=P$ (or more
precisely, a tubular neighborhood of $C$ in $P$).

Another obvious example is when $C$ is just a point $p$: then the
construction of Thm. \ref{emb} delivers as $\tP$ any slice
through $x$ transversal to the symplectic leaf $\cO_p$.

Now if $C_1\subset P_1$ and $C_2\subset P_2$ are pre-Poisson
submanifolds of Poisson manifolds, the cartesian product
$C_1\times C_2 \subset P_1\times P_2$ also is, and if the
construction of Thm. \ref{emb} gives cosymplectic submanifolds
$\tP_1\subset P_1$ and $\tP_2\subset P_2$, the same construction
applied to $C_1\times C_2$ (upon suitable choices of complementary
subbundles) delivers the cosymplectic submanifold $\tP_1\times
\tP_2$ of $P_1\times P_2$.
In particular, if $C_1$ is coisotropic and $C_2$ just a point $p$,
then $C_1\times \{p\}$ is pre-Poisson.
\end{ep}

The \emph{sufficient} condition above is not necessary (i.e. the
first inclusion in the diagram above is strict), as either of the
following simple examples shows.
\begin{ep}
Take $C$ to be the vertical line $\{x=y=0\}$ in the Poisson
manifold $(P,\Pi)=(\RR^3,f(z)\partial_x\wedge \partial_y)$, where
$f$ is any function with at least one zero. Then $C$ is a
Poisson-Dirac submanifold (with zero as induced Poisson structure),
hence taking $\tP:=C$ we obtain a Poisson-Dirac submanifold in
which $C$ embeds coisotropically. The sufficient conditions here
is not satisfied, for the rank of $TC + \sharp N^*C$ at $(0,0,z)$
is 3 at points where $f$ does not vanish and 1 at points where
$f$ vanishes.
\end{ep}

\begin{ep}
Consider the Poisson manifold $(P,\Pi)=(\RR^4,x^2\partial_x\wedge
\partial_y+z\partial_z\wedge \partial_w)$ as in Example 6
of \cite{CrF} and  the submanifold $C=\{(z^2,0,z,0):z\in \RR\}$.
The rank of $TC + \sharp N^*C$ is 3 away from the origin (because
there $C$ is an isotropic submanifold in an open symplectic leaf
of $P$) and 1 at the origin (since $\Pi$ vanishes there). The
submanifold $\tP=\{(z^2,0,z,w):z,w\in \RR \} $ is Poisson-Dirac
and it clearly contains $C$ as a coisotropic submanifold.
\end{ep}

The \emph{necessary} condition above is not  sufficient (i.e. the
second inclusion in the diagram above is strict):
\begin{ep}\label{mariusrui}
 In
Example 3 in Section 8.2 of \cite{CrF} the authors consider the
manifold $P=\CC^3$ with complex coordinates $x,y,z$. They specify a
Poisson structure on it by declaring the symplectic leaves to be
the complex lines given by $dy=0,dz-ydx=0$, the symplectic forms
being the restrictions of the canonical symplectic form on
$\CC^3$. They consider as submanifold $C$ the complex plane $\{z=0\}$
and show that $C$ is point-wise Poisson-Dirac (i.e. $TC\cap \sharp
N^*C=\{0\}$ at every point), but that the induced bivector field
is not smooth. Being point-wise Poisson-Dirac, $C$ satisfies the
necessary condition above. However there exists no Poisson-Dirac
submanifold $\tP$ of $P$ in which $C$ embeds coisotropically.
Indeed at points $p$ of $C$ where $y\neq 0$ we have $T_pC\oplus
T_p\cO=TP$ (where as usual $\cO$ is a symplectic leaf of $P$
through $p$), from which follows that $\sh|_{N_p^*C}$ is injective
and $T_pC\oplus \sh N_p^*C=TP$. From Lemma \ref{poisLA} (notice
that the subspace $R$ there must have trivial intersection with
$T_pC\oplus \sh N_p^*C$, so $R$ must be the zero subbundle over
$C$) it follows that the only candidate for $\tP$ is $C$ itself.
However, as we have seen, the Poisson bivector induced on $C$ is
not smooth. (More generally, examples are provided by any
submanifold $C$ of a Poisson manifold $P$ which is point-wise
Poisson-Dirac but not Poisson-Dirac and for
which there exists a point $p$ at which $T_pC\oplus T_p\cO=TP$.) \\
Notice that this provides an  example for the claim made in Remark
\ref{nons}, because the zero subbundle $R$ over $C$ satisfies equation \eqref{condPois} at every point of $C$ and is
obviously a smooth subbundle.
\end{ep}

We end with two examples of submanifolds $C$ which do not satisfy the necessary condition above. In particularly they can
not be imbedded coisotropically in any Poisson-Dirac submanifold.
\begin{ep} The submanifold $C=\{(x_1,x_2,x_2^2,x_1^2)\}$ of the
symplectic manifold $(P,\omega)=(\RR^4,dx_1\wedge dx_3+dx_2\wedge
dx_4)$ has characteristic distribution of rank 2 on the points
with $x_1=x_2$ and rank zero on the rest of $C$. The rank of the
characteristic distribution locally decreases, hence $C$ does not
satisfies the necessary condition above.
\end{ep}

\begin{remark}
If $C$ is a submanifold of a symplectic manifold $(P,\omega)$,
then the necessary and the sufficient conditions coincide, both
being equivalent to saying that the characteristic distribution of
$C$ (which can be described as $\ker(i^*_C\omega)$ for $i_C$ the
inclusion) have constant rank, i.e. that $C$ be presymplectic.
\end{remark}

\begin{ep}\label{noncont}
Consider the Poisson 
 manifold $(\RR^6,x_1\partial_{x_2}\wedge \partial_{x_4}+(\partial_{x_3}+x_1\partial_{x_5})\wedge\partial_{x_6})$. 
Let $C$ be the three-dimensional subspace given by setting $x_4=x_5=x_6=0$. The characteristic subspaces are all one-dimensional, spanned by $\partial_{x_3}$ at points of $C$ where $x_1=0$ and by $\partial_{x_2}$ on the rest of $C$. Hence the characteristic subspaces do not form a smooth distribution,
and can not be the image of the anchor map of any Lie algebroid over $C$. Therefore $C$ does not
satisfies the necessary condition above.
\end{ep}

\section{Reduction of submanifolds and deformation quantization of pre-Poisson submanifolds}\label{red}

In this section we consider the set of basic functions on a submanifold of a Poisson manifold,
and show that in certain cases it is a Poisson algebra and that it can be deformation quantized.

Given any submanifold $C$ of a Poisson manifold $(P,\Pi)$, it is
natural to consider the characteristic ``distribution'' $\sharp
N^*C \cap TC$ (which by Remark \ref{intr} consists of the kernels
of the restrictions to $C$ of the symplectic forms on the
symplectic leaves of $P$) and the set of basic functions on $C$ 
$$C^{\infty}_{bas}(C)=\{f  \in C^{\infty}(C): df|_{\sharp N^*C \cap
TC}=0\}.$$ 
 $\sharp
N^*C \cap TC$ usually does not have a constant rank and may not be smooth;
if it is and the quotient $\underline{C}$ is a smooth
manifold, then $C^{\infty}_{bas}(C)$ consists exactly of pullbacks of functions on
$\underline{C}$.

Let us endow $C$ with the (possibly non-smooth)
point-wise Dirac structure $i^*L_{P}$, where $i\colon C\rightarrow
P$ is the inclusion and $L_{P}$ is the Dirac structure
corresponding to $\Pi$. Then, 
since $\sharp N^*C \cap TC=i^*L_P\cap TC$,
$C^{\infty}_{bas}(C)$ is exactly the
set of basic functions of $(C,i^*L_P)$ in the sense of Dirac geometry.
 Given basic functions $f,g$ the expression
$$\{f,g\}_C(p):=Y(g)$$ is well-defined. Here $Y$ is any element of $T_pC$ such that
$(Y,df_p)\in i^*L_{P}$, and it 
  exists because
the annihilator of $i^*L_P\cap TC$ is the projection onto $T^*C$ of $i^*L_P$.
Notice that $C^{\infty}_{bas}(C)$ and $\{\bullet,\bullet\}_C$ are \emph{intrinsic} to $C$ in the following sense:
they depend only on the point-wise Dirac structure $i^*L_{P}$ on $C$, and
if $\bar{P}$ is a submanifold
of $(P,\Pi)$ containing $C$, $L_{\bar{P}}$ the point-wise Dirac structure on $\bar{P}$ induced by $P$ and
$\bar{i}:C \rightarrow \bar{P}$ the inclusion, then
$\bar{i}^*L_{\bar{P}}=i^*L_{P}$ by the functoriality of pullback. 

The expression
$\{f,g\}_C(p)$ does not usually vary
smoothly with $p$,
 so we can \emph{not} conclude that $C^{\infty}_{bas}(C)$ with $\{\bullet,\bullet\}_C$
is a Poisson algebra.
There is however a Poisson algebra that $C$ inherits from $P$
 \cite{CaFa1}, namely ${\cF} /(\cF\cap \cI)$, where
 $\cI$ denotes the set of functions on $P$ that vanish on
$C$ and $\cF:=\{\hat{f}\in C^{\infty}(P):\{\hat{f},\cI\}\subset
\cI\}$ (the so-called first class functions).
${\cF} /(\cF\cap \cI)$ is exactly the subset of functions $f$ on
$C$ which admits an extension to some function $\hat{f}$ on $P$
whose differential annihilates $\sharp N^*C$ (or equivalently
$X_{\hat{f}}|_C\subset TC$). The bracket of ${\cF} /(\cF\cap \cI)$ is computed as follows:
$$\{f,g\}=\{\hat{f},\hat{g}\}_P|_C=X_{\hat{f}}(g)|_C$$
for extensions as above. Notice that ${\cF} /(\cF\cap \cI)\subset C^{\infty}_{bas}(C)$,
and that the Poisson
bracket $\{\bullet,\bullet\}$ on ${\cF} /(\cF\cap \cI)$  coincides with  $\{\bullet,\bullet\}_C$:
if $f,g$ belong to ${\cF} /(\cF\cap \cI)$ we can compute
$\{f,g\}_C$ by choosing
 $Y=X_{\hat{f}}$ for some extension $\hat{f}\in \cF$.
 
 \begin{prop}\label{basic}
Let $C$ be any submanifold of a Poisson manifold $(P,\Pi)$. If
there exists  a Poisson-Dirac submanifold $\tP$ of $P$ in which
$C$ is contained coisotropically, then the set of basic functions
on $C$ has an intrinsic Poisson algebra structure, and
$({\cF} /(\cF\cap \cI),\{\bullet,\bullet\})$ is a Poisson subalgebra.
\end{prop}
\begin{proof}
We add a tilde  in the notation introduced above when we view $C$
as a submanifold of the Poisson manifold $\tP$ instead of $P$. By the last paragraph before the statement of this proposition, since
$\tilde{\sharp} N^*C\subset TC$,  it follows that
$\tilde{\cF}/\tilde{\cI}=C^{\infty}_{bas}(C)$.
So $(C^{\infty}_{bas}(C),\{\bullet,\bullet\}_C)$ is a Poisson algebra structure intrinsically associated to $C$, and it 
 contains ${\cF} /(\cF\cap \cI)$ as a Poisson subalgebra.\end{proof}

By Thm. \ref{emb} pre-Poisson submanifolds $C$ satisfy the
assumption of Prop. \ref{basic}, hence
  they admit a Poisson algebra structure on their
space of basic functions. This fact  was already established in  Theorem 3 of \cite{CaFa1},
where furthermore it is shown that ${\cF} /(\cF\cap \cI)$ is the whole space of basic functions.
 Now we state our result about deformation
quantization:

\begin{thm}\label{dqsub}
Let $C$ be a pre-Poisson submanifold, and assume that the first
and second Lie algebroid cohomology of $N^*C\cap \sharp^{-1}TC$
vanish. Then the Poisson algebra $C^{\infty}_{bas}(C)$, endowed with the bracket inherited from $P$, admits a
 deformation quantization.
\end{thm}
\begin{proof}
By Thm. \ref{emb} we can embed $C$ coisotropically in some
cosymplectic submanifold $\tP$. We invoke
Corollary 3.3 of \cite{CaFeCo2}: if the first and second Lie
algebroid cohomology of the conormal bundle of a coisotropic
submanifold vanish, then the Poisson algebra of basic functions on
the coisotropic submanifold (with the bracket inherited from the ambient  Poisson manifold, which
in our case is  $\tP$) admits a deformation quantization.
Now by Prop. \ref{basic} the
Poisson bracket  on $C^{\infty}_{bas}(C)$ induced by $P$ agrees with the one induced by the
embedding  in  $\tP$. 
Further the conditions in Corollary 3.3 of \cite{CaFeCo2} translate into the conditions stated in the proposition
because
 the conormal bundle of  $C$ in $\tP$ is
isomorphic to $N^*C\cap \sharp^{-1}TC$ as a Lie algebroid, see Prop.
\ref{algoid}.
\end{proof}

\section{Subgroupoids associated to pre-Poisson submanifolds}\label{groids}

Let $C$ be a pre-Poisson submanifold of a Poisson manifold
$(P,\Pi)$. In Prop. \ref{algoid} we showed that $N^*C\cap
\sharp^{-1}TC$ is a Lie subalgebroid of $T^*P$. When $\sharp N^*C$
has constant rank there is another Lie subalgebroid
associated  to $C$; it is obtained by taking the pre-image of $TC$ under the anchor map, i.e. it is  $\sharp^{-1} TC=(\sharp N^*C)^{\circ}$. Now we
  assume that $T^*P$ is an integrable Lie algebroid, i.e. that
 the source simply connected (s.s.c.) symplectic groupoid $(\Gs(P),\Omega)$ of $(P,\Pi)$
 exists. In this section we study the (in general only immersed) subgroupoids of
 $\Gs(P)$ integrating $N^*C\cap
\sharp^{-1}TC$ and $\sharp^{-1} TC$.  Here, for any Lie
subalgebroid $A$ of $T^*P$ integrating to a s.s.c. Lie groupoid
$G$, we take ``subgroupoid'' to mean the (usually just immersed)
image of the (usually not injective) morphism $G\rightarrow
\Gamma_s(P)$ induced from the inclusion $A\rightarrow T^*P$.

By Thm. \ref{emb} we can find a cosymplectic submanifold $\tP$ in
which $C$ lies coisotropically. We first make a few remarks on the
subgroupoid corresponding to $\tP$.

\begin{lemma}\label{ping}
The  subgroupoid of $\Gs(P)$ integrating $\sharp^{-1}T \tP$ is
$\stpttp$ and is a symplectic subgroupoid. Its source (target) map
is a Poisson (anti-Poisson) map onto $\tP$, where the latter is
endowed with the Poisson structure induced by $(P,\Pi)$.
\end{lemma}
\begin{proof}
  According to Thm. 3.7 of
\cite{Xu} the subgroupoid\footnote{In \cite{Xu} this is claimed
only when the subgroupoid integrating $(\sharp N^*\tP)^{\circ}$ is
an embedded subgroupoid, however the proof there is valid for
immersed subgroupoids too.} of $\Gs(P)$ corresponding to $\tP$,
i.e. the one integrating $(\sharp N^*\tP)^{\circ}$, is a
symplectic subgroupoid of $\Gs(P)$. It is given by $\stpttp$,
because  $(\sharp
N^*\tP)^{\circ}=\sharp^{-1}T\tP$.

To show that the maps $\stpttp \rightarrow \tP$ given
by the source and target maps of $\stpttp$ are Poisson
(anti-Poisson) maps proceed as follows. Take a function $\tilde{f}$ on $\tP$,
and extend it to a function $f$ on $P$ so that $df$ annihilates $\sharp N^*\tP$, i.e. so that
$X_{f}$ is tangent to
$\tP$ along $\tP$. 
Since $s:\Gs(P)\rightarrow P$ is a Poisson map and $s$-fibers are
symplectic orthogonal to $t$-fibers we know that the vector field
$X_{s^*f}$ on $\Gs(P)$ is tangent to $\stpttp$,
i.e. that $d(\bs^*f)$ annihilates $T(\stpttp)^{\Omega}$. Hence, denoting by
$\tilde{s}$ the source map of $\stpttp$, we have
$$\tilde{s}^*\{\tilde{f}_1,\tilde{f}_2\}=\tilde{s}^*(\{{f}_1,{f}_2\}|_{\tP})=
\{s^*f_1,s^*f_2\}|_{\stpttp}= \{\tilde{s}^*f_1,\tilde{s}^*f_1\},$$
i.e. $\tilde{s}$ is a Poisson map. A similar reasoning holds for
$\tilde{t}$.
\end{proof}

Now we describe the subgroupoid integrating $N^*C\cap
\sharp^{-1}TC$:
\begin{prop}\label{groid}
Let $C$ be a pre-Poisson submanifold of $(P,\Pi)$. Then the
subgroupoid of $\Gs(P)$ integrating $N^*C\cap \sharp^{-1}TC$ is an
isotropic subgroupoid of $\Gs(P)$.
\end{prop}
\begin{proof}
The canonical vector bundle isomorphism $i:T^*\tP\cong (\sharp
N^*\tP)^{\circ}$ is a Lie algebroid isomorphism, where $T^*\tP$ is
endowed with the cotangent algebroid structure coming from the
Poisson structure on $\tP$ (Cor.
2.11 and Thm. 2.3 of \cite{Xu}).
 Integrating this
 algebroid isomorphism we obtain a Lie groupoid morphism from $\Gs(\tP)$, the s.s.c.
Lie groupoid integrating $T^*\tP$, to $\Gs(P)$, and the image of
this morphism is
 $\stpttp$. Since by Lemma \ref{ping}
the symplectic form on $\stpttp$ is multiplicative, symplectic and  the source map is
a Poisson map,
pulling back the symplectic form on $\stpttp$
 endows $\Gs(\tP)$ with the structure
of the s.s.c. symplectic groupoid of $\tP$. The subgroupoid of
$\Gs(\tP)$ integrating $N^*_{\tP}C$, the annihilator of $C$ in
$\tP$, is Lagrangian (\cite{Ca}, Prop. 5.5). Hence
$i(N^*_{\tP}C)$, which by Prop. \ref{algoid} is equal to $N^*C\cap
\sharp^{-1}TC$, integrates to a Lagrangian subgroupoid of
$\stpttp$, which therefore is an isotropic subgroupoid of
$\Gs(P)$.
\end{proof}

Now we consider $\sharp^{-1}TC$. For any submanifold $N$,
$\sharp^{-1}TN$ has constant rank iff it is a Lie subalgebroid of
$T^*P$, integrating to the subgroupoid $\bs^{-1}(N)\cap
\bt^{-1}(N)$ of $\Gs(P)$. So the constant rank condition on
$\sharp^{-1}TN$ corresponds to a smoothness condition on
$\bs^{-1}(N)\cap \bt^{-1}(N)$.

\begin{remark}\label{onetwo} 1) If $\sharp^{-1}TN$ has constant rank it follows
that the Poisson structure on $P$ pulls back to a smooth Dirac
structure on $N$, and that $\bs^{-1}(N)\cap \bt^{-1}(N)$ is an
over-pre-symplectic
 groupoid inducing the same Dirac
structure on $N$ (Ex. 6.7 of \cite{BCWZ}). Recall from Def. 4.6 of \cite{BCWZ}
that an over-pre-symplectic groupoid is a Lie groupoid $G$ over a base
$M$ equipped with a closed multiplicative 2-form $\omega$ such
that $ker \omega_x\cap ker (d\bs)_x  \cap ker (d\bt)_x $ has rank $dim
G-2 dim M$ at all $x\in M$. Further, $\bs^{-1}(N)\cap
\bt^{-1}(N)$ has dimension equal to $2dim N+rk(N^*N\cap N^*\cO)$,
where $\cO$ denotes any symplectic leaf of $P$ intersecting $N$.

\noindent 2) For a pre-Poisson submanifold $C$, the condition that
$\sharp^{-1}TC$ have constant rank is equivalent to the
characteristic distribution $TC\cap \sh N^*C$ having constant
rank.  This follows trivially from
$rk (\sharp N^*C+TC)= rk (\sh N^*C)+ \dim C - rk (TC\cap \sh N^*C)
$.
\end{remark}

\begin{prop}\label{st}
Let $C$ be a pre-Poisson submanifold with constant-rank
characteristic distribution. Then for any cosymplectic submanifold
$\tP$ in which $C$ embeds coisotropically, $\sctc$ is a
coisotropic subgroupoid of $\stpttp$.
\end{prop}
\begin{proof}
By the comments above we know that $\sharp^{-1}TC$ is a Lie
subalgebroid, hence $\sctc$ is a (smooth) subgroupoid of $\Gs(P)$,
and it is clearly contained in $\stpttp$.
 We saw in Lemma \ref{ping}
 that $\stpttp$ is endowed with a symplectic form for which
  its source and target maps are (anti-)Poisson
maps onto $\tP$. Further its source and target fibers symplectic
orthogonals of each other. Since $C\subset \tP$ is coisotropic,
this implies
that $\sctc$ is coisotropic in $\stpttp$.
\end{proof}

We now describe the subgroupoids corresponding to pre-Poisson
manifolds.

\begin{prop}\label{pre}
Let $C$ be any submanifold of $P$.  Then $\sctc$ is a (immersed)
presymplectic submanifold of $\Gs(P)$ iff $C$ is pre-Poisson and its
characteristic distribution has constant rank. In this case the
characteristic distribution of $\sctc$ has rank $2rk(\sharp N^*C
\cap TC)+rk(N^*C\cap N^*\cO)$, where $\cO$ denotes the symplectic
leaves of $P$ intersecting $C$.
\end{prop}

\begin{proof}
Assume that $\sctc$ is a (immersed) presymplectic submanifold of $\Gs(P)$.
  We apply the same proof as in Prop. 8 of \cite{CrF}:
there is an isomorphism of vector bundles $T\Gs(P)|_P\cong
TP\oplus T^*P$, under which the non-degenerate bilinear form
$\Omega|_P$ corresponds to $ (X_1\oplus \xi_1, X_2\oplus \xi_2)
:=\langle \xi_1,X_2 \rangle - \langle \xi_2,X_1 \rangle +
\Pi(\xi_1,\xi_2)$. Under the above isomorphism $T(\sctc)$
corresponds to $TC\oplus \sharp^{-1}TC$, and a short computation
shows that the restriction of $(\bullet,\bullet)$ to $TC\oplus
\sharp^{-1}TC$ has kernel $(TC\cap \sh N^*C)\oplus
(\sharp^{-1}TC\cap N^*C)$, which therefore has constant rank. From the smoothness of  $\sctc$ it
follows that $(\sh N^*C)^{\circ}=\sharp^{-1}TC$ has constant rank.
This has two consequences: first by Remark \ref{onetwo}
$C$ has
characteristic distribution of constant rank. Second, the above
  kernel is a direct sum of two intersections of smooth
subbundles, so $\sharp^{-1}TC\cap N^*C$ has constant rank,
i.e. (taking annihilators) $C$ is pre-Poisson.
 
The other direction follows from Prop. \ref{st}.
\end{proof}

\begin{remark}One can wonder whether any subgroupoid of a symplectic groupoid
$(\Gs(P),\Omega)$ which is a presymplectic submanifold (i.e.
$\Omega$ pulls back to a constant rank 2-form) is contained
coisotropically in some symplectic subgroupoid of $\Gs(P)$. This
would be exactly the ``groupoid'' version of Thm. \ref{emb}. The
above Prop. \ref{st}  and Prop. \ref{pre} together tell us that
this is the case when the subgroupoid has the form $\sctc$, where
$C\subset P$ is its base. In general the answer to the above
question is negative, as the following counterexample shows.

Let $(P,\omega)$ be some simply connected symplectic manifold, so
that $\Gs(P)=(P\times P,\omega_1-\omega_2)$ and the units are
embedded diagonally. Take  $C$ to be any 1-dimensional closed
submanifold of $P$. $C {\rightrightarrows} C$ is clearly  a
subgroupoid and a presymplectic submanifold; since
$\omega_1-\omega_2$ there pulls back to zero, any subgroupoid $G$
of $P\times P$ in which $C {\rightrightarrows} C$ embeds
coisotropically must have dimension 2. If the base of $G$ has
dimension 2 then $G$ is contained in the identity section of
$P\times P$, which is Lagrangian. So let us assume that the base
of $G$ is $C$. Then $G$ must be contained in $C \times C$, on
which $\omega_1-\omega_2$ vanishes because $C\subset P$ is
isotropic. So we conclude that there is no symplectic subgroupoid
of $P\times P$ containing $C {\rightrightarrows} C$ as a
coisotropic submanifold.
\end{remark}

\section{Existence of coisotropic embeddings of Dirac manifolds in Poisson
manifolds}\label{emdpois}

 Let $(M,L)$ be a smooth Dirac
manifold. We ask when $(M,L)$ can be embedded coisotropically in
some Poisson manifold $(P,{\Pi})$, i.e. when there exists an
embedding $i$ such that $i^* L_{{P}}=L$ and $i(M)$ is a
coisotropic submanifold of $P$. Notice that for arbitrary coisotropic
embeddings $i^* L_{{P}}$ is usually not even continuous: for example the $x$-axis in $(\RR^2,x\partial_x\wedge\partial_y)$ is coisotropic,
but the pullback  structure is not continuous at the origin.

 When $M$ consists of exactly one leaf, i.e.
when $M$ is a manifold endowed with a closed 2-form $\omega$, the
existence and uniqueness of
coisotropic embeddings in symplectic manifolds was considered by
Gotay in the short paper \cite{Go}: the coisotropic embedding exists
if{f} $\ker \omega$ has
constant rank, and in that case one has uniqueness up to neighborhood equivalence.
 Our strategy will be to check
if we can apply Gotay's arguments ``leaf by leaf'' \emph{smoothly}
over $M$. Recall that
$L\cap TM$ is the kernel of the 2-forms on the presymplectic
leaves of $(M,L)$.

\begin{thm}\label{E}
$(M,L)$ can be embedded coisotropically in a Poisson manifold
if{f} $L\cap TM$ has constant rank.
\end{thm}
\begin{proof}
Suppose that an embedding $ M\rightarrow P$ as above
exists. Then $L\cap TM$ is equal ${\sharp} N^*C$ (where $N^*C$ is
the conormal bundle of $C$ in $P$), the image of a vector bundle under
a smooth bundle map, hence its rank can locally only increase. On
the other hand the rank of
 $L\cap TM$, which is the intersection of two
 smooth bundles, can locally only decrease.
 Hence the rank of $L\cap TM$ must be constant on $M$.

 Conversely, assume that the rank of $E:=L\cap TM$ is constant
and define $P$ to be the total space of the vector bundle
$\pi\colon E^*\rightarrow M$. We define the Poisson structure on
$P$ as follows. First take the pullback Dirac structure $\pi^*L$
(which is smooth and integrable since $\pi$ is a submersion). Then
 choose a smooth
distribution $V$ such that $E\oplus V= TM$. This choice gives an
embedding $i\colon E^* \rightarrow T^*M$, which we can use to
pull back the canonical symplectic form $\omega_{T^*M}$. Our
Poisson structure is $L_{E^*}:=\tau_{i^* \omega_{T^*M}}\pi^*L$, i.e. it is
obtained applying to $\pi^*L$ the gauge transformation\footnote{
Given a Dirac structure $L$ on a vector space $W$, the gauge transformation of $L$ by a bilinear form $B\in \wedge^2 W^*$ is   $\tau_BL:=\{(X,\xi + i_XB):(X,\xi)\in L\}$.
Given a Dirac structure $L$ on a manifold, the gauge-transformation $\tau_BL$ by closed 2-form $B$
is again a Dirac structure (i.e. $\tau_BL$ is again closed under the Courant bracket).}
  by the closed 2-form $i^*
\omega_{T^*M}$. It is clear that $L_{E^*}$ is a smooth Dirac structure; we still have to show that
it is actually Poisson, and that the zero section is coisotropic.
In more concrete terms $(E^*,L_{E^*})$ can be described as
follows: the leaves are all of the form $\pi^{-1}
(F_{\alpha})$ for $(F_{\alpha},\omega_{\alpha})$ a presymplectic
leaf of $M$. The 2-form on the leaf is given by adding to
$(\pi|_{\pi^{-1}(F_{\alpha})})^*\omega_{\alpha}$ the 2-form
$i_{\alpha}^*\omega_{T^* F_{\alpha}}$. The latter is defined
considering  the distribution $V\cap TF_{\alpha}$  transverse  to
$E|_{F_{\alpha}}$ in $TF_{\alpha}$, the induced embedding
$i_{\alpha}\colon \pi^{-1}(F_{\alpha})=E^*|_{F_{\alpha}} \rightarrow T^*
F_{\alpha}$, and pulling back the canonical symplectic form.
(One can check that $i_{\alpha}^*\omega_{T^* F_{\alpha}}$ is the pullback of
$i^* \omega_{T^*M}$ via the inclusion of the leaf in $E^*$).
But
this is exactly Gotay's recipe to endow (an open subset of) $\pi^{-1} (F_{\alpha})$
with a symplectic form so that $F_{\alpha}$ is embedded as a
coisotropic submanifold. Hence we conclude that a neighborhood of the zero section of $E^*$, with the
above Dirac structure, is actually a Poisson manifold and that $M$
is embedded as a coisotropic submanifold.
\end{proof}

We comment on how choices affect the construction of Thm. \ref{E}.
We need the following version of Moser's theorem for Poisson
structures (see Section 3.3. of \cite{AM}) : suppose we are given
 Poisson structures $\Pi_t$ on some manifold $P$, $t\in [0,1]$. Assume
 that each $\Pi_t$ is related to $\Pi_0$ via the gauge transformation
by some closed 2-form $B_t$, i.e. $\Pi_t=\tau_{B_t}\Pi_0$. This
means that  the symplectic foliations agree and on each symplectic
leaf $\cO$ we have $\Omega_t=\Omega_0+i_{\cO}^*B_t$, where
$\Omega_0$,$\Omega_t$ are the symplectic forms on the leaf $\cO$
and $i_{\cO}$ the inclusion. Assume further that each
$\frac{d}{dt} B_t$ be exact, and let $\alpha_t$ be a smooth family
of primitives vanishing on some submanifold $M$. Then the time-1
flow of the Moser vector field
 $\sharp_t \alpha_t$ is
defined in a tubular neighborhood of $M$, it fixes $M$ and maps
$\Pi_0$ to $\Pi_1$. Here  $\sharp_t$ denotes
the map $T^*P\rightarrow TP$ induced by $\Pi_t$.

\begin{prop}\label{splitting}
 Different choices of splitting
$V$ in the construction of Thm. \ref{E} yield  
isomorphic Poisson structures on  $E^*$. Hence, given a Dirac
manifold $(M,L)$ for which $L\cap TM$ has constant rank, there is
a canonical (up to neighborhood equivalence) Poisson manifold in
which $M$ embeds coisotropically.
\end{prop}
\begin{proof}
Let $V_0,V_1$ be two different splittings as in Thm. \ref{E}, i.e.
$E\oplus V_i=TM$ for $i=0,1$. We can interpolate between them by
defining the graphs $V_t:=\{v+tAv:v\in V_0\}$ for $t\in [0,1]$,
where $A \colon V_0\rightarrow E$ is determined by requiring that
its graph be $V_1$. Obviously each $V_t$  also gives a  splitting
$E\oplus V_t=TM$; denote by $i_t\colon E^* \rightarrow T^*M$ the
corresponding embedding. We obtain Dirac structures $\tau_{i_t^*
\omega_{T^*M}}\pi^*L$ on the total space of $\pi \colon E^*
\rightarrow M$; by Thm. \ref{E} they correspond to Poisson
bivectors, which we denote by $\Pi_t$. These Poisson structures
are related by a  gauge transformation: $\Pi_t=\tau_{B_t} \Pi_0$ 
for $B_t:=i_t^*
\omega_{T^*M}-i_0^* \omega_{T^*M}$. A primitive of $\frac{d}{dt}B_t$ is given by $\frac{d}{dt}i_t^*
\alpha_{T^*M}$; notice that this primitive vanishes at points of
$M$, because the canonical 1-form $\alpha_{T^*M}$ on $T^*M$
vanishes along the zero section. Hence the time-1 flow of $\sharp_t
(\frac{d}{dt}i_t^* \alpha_{T^*M})$ fixes $M$ and maps  $\Pi_0$ to
$\Pi_1$.
\end{proof}

Assuming that $(M,L)$ is integrable we describe the   symplectic
groupoid of $(E^*,L_{E^*})$, the Poisson
manifold constructed in Thm. \ref{E} with a choice of
distribution $V$. It is $\pi^*(\Gs(M))$, the pullback via
$\pi\colon E^*\rightarrow M$ of the presymplectic groupoid of $M$,
endowed with the following symplectic form: the pullback via
$\pi^*(\Gs(M))\rightarrow \Gs(M)$ of the presymplectic form on the
groupoid $\Gs(M)$, plus
$\bs^*(i^*\omega_{T^*M})-\bt^*(i^*\omega_{T^*M})$, where
$i\colon E^* \rightarrow T^*M$ is the inclusion given by the
choice of distribution $V$, $\omega_{T^*M}$ is the canonical
symplectic form, and $\bs,\bt$ are the source and target maps of
$\pi^*(\Gs(M))$. This follows easily from Examples 6.3 and 6.6 in
\cite{BCWZ}. Notice that this groupoid is source simply connected
when $\pi^*(\Gs(M))$ is.

Now we can give
an affirmative answer to the possibility raised in \cite{CrF}
(Remark (e) in Section 8.2), although we prove it ``working
backwards``; this is the ``groupoid'' version of Gotay's embedding theorem.
Recall  that a presymplectic groupoid is a Lie groupoid
$G$ over $M$ with $dim(G)=2dim(M)$ equipped with a closed multiplicative 2-form $\omega$
such that $ker \omega_x\cap ker (d\bs)_x  \cap ker (d\bt)_x =0 $   at
all $x\in M$  (Def.
2.1 of \cite{BCWZ}).

\begin{prop}
Any presymplectic groupoid   with constant rank characteristic
distribution can be embedded coisotropically as a Lie subgroupoid
in a symplectic groupoid.
\end{prop}

\begin{proof}
By  Cor. 4.8 iv),v) of
\cite{BCWZ}, a presymplectic groupoid $\Gs(M)$ has characteristic
distribution (the kernel of the multiplicative 2-form) of constant
rank if{f} the Dirac structure $L$ induced on its base $M$ does.
We can embed $(M,L)$
coisotropically in the Poisson manifold $(E^*,L_{E^*})$
constructed in Thm. \ref{E}; we just showed that $\pi^*(\Gs(M))$ is a symplectic groupoid
for $E^*$.
$\Gs(M)$ embeds in $\pi^*(\Gs(M))$ as $\bs^{-1}(M)\cap
\bt^{-1}(M)$, and this embedding preserves both the groupoid structures and
the 2-forms.  $\bs^{-1}(M)\cap
\bt^{-1}(M)$ is a coisotropic subgroupoid of
$\pi^*(\Gs(M))$ because $M$ lies coisotropically in $E^*$ and
$\bs,\bt$ are (anti)Poisson maps.
\end{proof}

\begin{remark}
A partial converse to this proposition is given as follows: if
$\bs^{-1}(M)\cap \bt^{-1}(M)$
 is a coisotropic  subgroupoid of a symplectic
groupoid $\Gs(P)$, then $M$  is a coisotropic submanifold of the
Poisson manifold $P$, it has an smooth Dirac structure (induced
from $P$) with characteristic distribution of constant rank, and
$\bs^{-1}(M)\cap \bt^{-1}(M)$ is a \emph{over}-pre-symplectic
groupoid over $M$ inducing the same Dirac structure. This follows
from our arguments in section \ref{groids}.
\end{remark}

Now we draw the conclusions about deformation quantization. Recall
that for any Dirac manifold $(M,L)$ the set of admissible
functions \begin{equation}\label{adm} C^{\infty}_{adm}(M)=\{f\in
C^{\infty}(M): \text{ there exists a smooth vector field }X_f\text{
s.t. }(X_f,df)\subset L\}\end{equation} is naturally a Poisson
algebra \cite{Cou}, with bracket $\{f,g\}_M=X_f(g)$.

\begin{thm}\label{dqdirac}
Let $(M,L)$ be a Dirac manifold such that $L\cap TM$ has constant
rank, and denote by $\cF$ the regular foliation integrating $L\cap
TM$. If the first and second foliated de~Rham cohomologies of the
foliation $\cF$ vanish then the Poisson algebra of admissible
functions on $(M,L)$ admits a deformation quantization.
\end{thm}
\begin{proof}
 By Thm. \ref{E} we can embed $(M,L)$
coisotropically in a Poisson manifold $P$; hence we can apply
  Corollary 3.3 of \cite{CaFeCo2}: if the first and second Lie
algebroid cohomology of the conormal bundle of a coisotropic
submanifold vanish, then the Poisson algebra of basic functions on
the coisotropic submanifold admits a deformation quantization.
Since $L\cap TM$ has constant rank, the image of $L$ under $TM\oplus T^*M \rightarrow T^*M$  has
constant rank., so
the inclusion
$C^{\infty}_{adm}(M)\subset C^{\infty}_{bas}(M)$ is an equality. Further the Poisson algebra structure $\{\bullet,\bullet\}_M$
 on
$C^{\infty}_{bas}(M)$ coming from $(M,L)$ coincides with the one
induced by $M$ as a coisotropic submanifold of $P$, as follows
from Prop. \ref{basic} and $i^* L_{{P}}=L$. So
 when the assumptions are satisfied we really
deformation quantize $(C^{\infty}_{adm}(M),\{\bullet,\bullet\}_M)$.

Notice that in Thm. \ref{E} we constructed a Poisson manifold $P$
 of minimal
dimension, i.e. of dimension $dim M+rk (L\cap TM)$. The anchor map
$\sharp$ of the Lie algebroid $N^*C$ is injective, hence the Lie algebroids $N^*C$
and $L\cap TM$ are isomorphic. This allows us to state the
assumptions of Corollary 3.3 of \cite{CaFeCo2} in terms of the
foliation $\cF$ on $M$.
\end{proof}

\begin{prop}\label{li}
Let $(M,L)$ be a Dirac manifold such that $L\cap TM$ has constant
rank, and denote by $\cF$ the regular foliation integrating $L\cap
TM$. Then the foliated de~Rham complex $\Omega^{\bullet}_{\cF}(M)$
admits the structure of an $L_{\infty}$-algebra\footnote{The
$\lambda_n$ are derivations w.r.t. the wedge product, so one
 actually obtains what in \cite{CaFeCo2} is called
 a $P_{\infty}$ algebra.} $\{\lambda_n\}_{n\ge 1}$, the
differential $\lambda_1$ being the foliated de~Rham differential
and the bracket $\lambda_2$ inducing on
$H^0_{\lambda_1}=C^{\infty}_{bas}(M)$ the natural bracket
$\{\bullet,\bullet\}_M$. This $L_{\infty}$ structure is canonical up to 
$L_{\infty}$-isomorphism.
\end{prop}
\begin{proof}
By the proof of Thm. \ref{E} we know that $M$ can be embedded
coisotropically in a Poisson manifold $P$ so that the Lie
algebroids $N^*M$ and $L\cap TM$ are isomorphic.  After choosing
an embedding of $NM:=TP|_M/TM$ in a tubular neighborhood
of $M$ in $P$, Thm. 2.2 of \cite{CaFeCo2} gives the desired
$L_{\infty}$-structure. By Prop. \ref{splitting} the Poisson manifold $P$ is canonical up to neighborhood equivalence, so
the $L_{\infty}$-structure depends only on the choice of embedding of $NM$ in $P$;
the first author and Sch\"atz showed in \cite{CS} that different embeddings give the same
structure up to $L_{\infty}$-isomorphism.\end{proof}

\section{Uniqueness of coisotropic embeddings of Dirac
manifolds}\label{aissa}

The coisotropic embedding
of
Gotay \cite{Go} is unique up to
neighborhood equivalence, i.e. any two coisotropic embeddings of
a  fixed presymplectic manifold in symplectic manifolds
are
intertwined by a symplectomorphism which is the identity on the
coisotropic submanifold. It is natural to ask whether, given a Dirac manifold $(M,L)$
such that $L\cap TM$ have constant rank,
 the coisotropic embedding constructed in Thm. \ref{E} is
the only one up to neighborhood equivalence. In general the answer
will be negative: for example the origin is a coisotropic submanifold in
$\RR^2$ endowed either with the zero Poisson structure or with the
Poisson structure $(x^2+y^2)\partial_x\wedge
\partial_y$, and the two Poisson structures are clearly not
equivalent.

 As Aissa Wade pointed out to us, it is necessary to
require that the  Poisson manifold in which we embed
be of minimal dimension, i.e. of dimension $\dim M+\emph{rk}(L\cap TM)$.
Before presenting some partial results on the uniqueness problem we need a simple lemma.  
\begin{lemma}\label{easylemma}
Let $M$ be a coisotropic subspace of a Poisson vector space $(P,\Pi)$. Then 
  $codim(M)=dim(\sharp M^{\circ})$ if{f} $\sharp|_{M^{\circ}}$ is an injective map
  if{f}
   $M$ intersects transversely $\cO:=\sharp P^*$.
\end{lemma}
\begin{proof} The first equivalence is obvious by dimension reasons.
 For the second one notice that   $\sharp|_{M^{\circ}}$ is injective  
iff $M^{\circ}\cap \cO^{\circ}=\{0\}$, which taking annihilators is exactly  the transversality statement.
\end{proof}

\subsection{Infinitesimal uniqueness and global issues}

We  apply the construction of Gotay's uniqueness  proof \cite{Go} on each presymplectic leaf
of the Dirac manifold $M$; then we    show that under certain assumptions the
resulting diffeomorphism varies smoothly from leaf to leaf.

We start establishing infinitesimal uniqueness, for which we need a Poisson linear algebra lemma.

\begin{lemma}\label{bundleLA}
Let $(P,\Pi)$ be a Poisson vector space and $M$ a coisotropic subspace for which 
$dim(\sharp M^{\circ})=codim(M)$. Let $V$ be a complement to $E:=\sharp M^{\circ}$ in $M$.

There exists an isomorphism of Poisson vector spaces 
$$P \cong V\oplus E\oplus E^*$$ 
fixing $M$,
where the Poisson structure on the r.h.s. is such that the induced symplectic vector space is $((V\cap \cO)\oplus (E\oplus E^*), \Omega|_{V\cap \cO}\oplus \omega_E)$. Here $(\cO,\Omega)$ is the symplectic subspace corresponding to $(P,\Pi)$ and $\omega_E$ is the antisymmetric pairing on $E\oplus E^*$.
\end{lemma}

\begin{proof}
We claim first that $V \oplus \sharp
V^{\circ}=P$: indeed $V\cap \cO$ is a
symplectic subspace of $(\cO,\Omega)$,
being transverse to $E=ker(\Omega|_{\cO\cap M})$.
 Hence $(V\cap 
\cO)^{\Omega}$, which by   section \ref{def} is equal to
$\sharp V^{\circ}$, is a complement to $V\cap \cO$ in
$\cO$, so $V \oplus \sharp V^{\circ}=V+ \cO$, which equals $P$ by Lemma \ref{easylemma}.

 Now
we mimic the construction of  Gotay's uniqueness proof \cite{Go}:
since $E$ is Lagrangian in the symplectic subspace $\sharp
V^{\circ}$, by choosing a complementary lagrangian we can find a linear symplectomorphism $(\sharp
V^{\circ}, \Omega|_{\sharp V^{\circ}}) \cong (E\oplus E^*,
\omega_E)$ which is the identity on $E$. Adding to this $Id_V$ we obtain an isomorphism 
$$P=V \oplus \sharp V^{\circ}
\cong V\oplus E\oplus E^*,$$ which preserves the Poisson bivectors
 because  it restricts to an isomorphism $\cO
\cong (V\cap \cO)\oplus (E\oplus E^*)$ which  matches the symplectic forms $\Omega$
and $\Omega|_{V\cap \cO}\oplus \omega_E$. 
\end{proof}

\begin{prop}\label{bundle}
Suppose we are given a  Dirac manifold $(M,L)$ for which $L\cap
TM$ has constant rank $k$, and let $(P_1,\Pi_1)$ and $(P_2,\Pi_2)$ be
Poisson manifolds of dimension $\dim M+ k$ in
which $(M,L)$ embeds coisotropically. Then there is an isomorphism of
Poisson vector bundles $\Phi \colon TP_1|_M \rightarrow TP_2|_M$
which is
the identity on $TM$.
\end{prop}
\begin{proof} We choose a smooth distribution $V$ on $M$ completementary to $E:=L\cap TM$. 
For $i=1,2$, at every point $x\in M$ we  apply the construction of Lemma \ref{bundleLA} to the coisotropic subspace $T_xC$ of $T_xP_i$, obtaining
 smooth isomorphisms of Poisson vector bundles  $$TP_1|_M\cong V\oplus E\oplus E^*\cong TP_2|_M.$$ 
  Notice that middle Poisson vector bundle depends only on $(M,L)$ and $V$,  because for any symplectic 
  leaf $(\cO,\Omega)$ of $P_1$ or $P_2$
the bilinear form $\Omega|_{V_x\cap T_x\cO}$ is determined by the presymplectic form on the
 presymplectic leaf $\cO\cap M$ of $(M,L)$.
 \end{proof}

Making a regularity assumption we  can extend the infinitesimal
uniqueness of Prop. \ref{bundle} to a global statement.
\begin{prop}\label{moser}
Let $M$,$P_1$ and $P_2$ be as in Proposition \ref{bundle}, and
assume additionally that the presymplectic leaves of $(M,L)$ have
constant dimension. Then $P_1$ and $P_2$ are neighborhood
equivalent.
\end{prop}

\begin{proof}
The symplectic leaves of each $P_i$  have constant
dimension in a tubular neighborhood of $P_i$,
because they are transverse to $M$ by Lemma
 \ref{easylemma} and because of the assumption on the  presymplectic leaves of $(M,L)$.  
By choosing normal
bundles $N_i\subset TP_i|_M$ tangent to the symplectic leaves of $P_i$
  we can find identifications $\phi_i$ between the normal
bundles $N_i$ and tubular neighborhoods of $M$ in
$P_i$ which, for
every presymplectic leaf $F$ of $M$, identify  ${N_i}|_{F}$ and  the corresponding symplectic leaf of $P_i$.

Using the Poisson vector bundle isomorphism $\Phi\colon TP_1|_M
\rightarrow  TP_2|_M$ of Proposition \ref{bundle} we obtain an
identification $\phi_2\circ\Phi\circ \phi_1^{-1}$ between tubular
neighborhoods of $M$ in $P_1$ and $P_2$. Using this identification
 can view $\Pi_2$ as a Poisson structure
on $P:=P_1$ with two properties: it induces exactly the same
foliation as $\Pi_1$, and it coincides with $\Pi_1$ on $TP|_M$. We
want to show that there is a diffeomorphism near $M$, fixing $M$,
which maps $\Pi_1$ to $\Pi_2$.

To this aim we  apply Moser's theorem on each symplectic
leaf $\cO$ of $P$ (Thm. 7.1 of \cite{Ana}), in a way that varies smoothly with $\cO$. Denote by $\Omega_i$ the
symplectic form given by $\Pi_i$ on a leaf $\cO$. The convex
linear combination $(1-t)\Omega_1 +t\Omega_2$ is symplectic
(because $\Omega_1$ and $\Omega_2$ coincide at points of $M$). Let $F:=M\cap \cO$,
 identify $N|_F$ with a neighborhood in $\cO$ via $\phi_1$, and consider the retraction $\rho_t\colon N|_{F}\rightarrow
N|_{F},v\mapsto tv$
 where $t\in
[0,1]$.
Let $Q$ be the homotopy operator given by the 
 retraction $\rho_t$ (see Chapter 6 of \cite{Ana}); it satisfies $dQ-Qd=\rho_1^*-\rho^*_0$.
  So
$\mu:=Q(\Omega_2-\Omega_1)$ is a primitive for
$\Omega_2-\Omega_1$; furthermore $Q$ can be chosen so that $\mu$ vanishes at points of $F$. Consider the Moser vector field, obtained
inverting via $(1-t)\Omega_1 +t\Omega_2$ the 1-form $\mu$.
Following from time 1 to time 0 the flow of the Moser vector field
gives  a diffeomorphism  $\psi$ of  $\cO$ fixing $F$ such that
$\psi^*\Omega_2=\Omega_1$.

Notice that since the symplectic foliation of $P$ is regular near $M$ this construction varies smoothly from leaf to leaf.  Hence we
obtain a diffeomorphism $\psi$ of a tubular neighborhood of $M$ in $P$,
fixing $M$, which maps $\Pi_1$ to $\Pi_2$.
\end{proof}

 Since local uniqueness holds (see subsection
\ref{uniloc}) and since by Proposition \ref{bundle} there is no topological obstruction,
it seems that the
global uniqueness statement of Prop. \ref{moser} should hold in full generality (i.e. without the assumption on the
presymplectic foliation of $(M,L)$)
however we are not able to prove this.

The argument from \cite{AM} just before our Prop. \ref{splitting}
shows that the uniqueness of (minimal dimensional) coisotropic
embeddings of a given Dirac manifold $(M,L)$ is equivalent to the
following: whenever $(P_1,\Pi_1)$ and $(P_2,\Pi_2)$ are minimal
Poisson manifolds in which $(M,L)$ embeds coisotropically there
exists a diffeomorphism $\phi\colon P_1\rightarrow P_2$ near $M$
so that $\Pi_2$ and $\phi_*\Pi_1$ differ by the gauge
transformation by a closed 2-form $B$ vanishing on $M$. One could
hope that if $\phi\colon P_1\rightarrow P_2$ is chosen to match
symplectic leaves and to match $\Pi_1|_M$ and $\Pi_2|_M$ then
 a 2-form
$B$ as above automatically exists. This is
 not the case, as the following example shows.

\begin{ep}\label{r4}
Take $M =\RR^3$ with Dirac structure $$L=span\{(-x_1^2
\partial_{x_2},dx_1),
(x_1^2 \partial_{x_1}, dx_2), (\partial_{x_3},0)\}.$$ There are
two open presymplectic leaves $(\RR_{\pm}\times \RR^2,
\frac{1}{x_1^2}dx_1\wedge dx_2)$ and 1-dimensional presymplectic
leaves $\{0\}\times \{c\}\times \RR$ with zero presymplectic form
(for every real number $c$); hence our Dirac structure is a
product of the Poisson structure $ x_1^2
\partial_{x_1}\wedge
\partial_{x_2}$ and of the zero presymplectic form on the $x_3$-axis.
The characteristic distribution $L\cap TM$ is always $span
\partial_{x_3}$. Clearly the construction of Thm. \ref{E} gives
$$P_1:= (\RR^4, x_1^2
\partial_{x_1}\wedge
\partial_{x_2}+\partial_{x_3}\wedge
\partial_{y_3})$$ where $y_3$ is the coordinate on the fibers of
$P_1\rightarrow M$.

Another Poisson structure on $\RR^4$ with the same foliation as
$\Pi_1$ and which coincides with $\Pi_1$ along $M$ is the
following:
$$\Pi_2:= x_1^2
\partial_{x_1}\wedge
\partial_{x_2}+
\partial_{x_3}\wedge
\partial_{y_3}+x_1y_3\partial_{x_2}\wedge
\partial_{x_3}.$$
On each of the two open symplectic leaves  $\RR_{\pm}\times \RR^3$
the symplectic form corresponding to $\Pi_1$ is
$\Omega_1=\frac{1}{x_1^2}dx_1\wedge dx_2+dx_3\wedge dy_3$, whereas
the one corresponding to $\Pi_2$ is
$\Omega_2=\Omega_1+\frac{y_3}{x_1}dx_1\wedge dy_3$. Clearly the
difference $\Omega_1-\Omega_2$ does not extend to smooth a 2-form on the
whole of $\RR^4$. Hence there is no smooth 2-form on $\RR^4$
relating $\Pi_1$ and $\Pi_2$.

Nevertheless $\Pi_1$ and $\Pi_2$ are Poisson diffeomorphic: an explicit Poisson diffeomorphism is given by
the
global coordinate change that
transforms $x_2$ into
$x_2+\frac{y_3^2}{2}x_1$ and leaves the other coordinates
untouched.
\end{ep}

\subsection{Local uniqueness}\label{uniloc}

While we are not able to prove a global uniqueness statement
in the general case, we prove in this
subsection that  local uniqueness  holds. We start with
a normal form statement.

\begin{prop}\label{can}
Let $M^m$ be a coisotopic submanifold of a Poisson manifold $P$ such that $k:=codim(M)$ equals $rk(\sharp N^*M)$. Then
about any $x\in M$ there is a neighborhood $U\subset P$ and
coordinates $\{q_1,\dots,q_k,p_1,\dots,p_k,y_1,\dots,y_{m-k}\}$
defined on $U$ such that locally $M$ is given by the constraints
$p_1=0,\dots,p_k=0$ and
\begin{equation}\label{canform} \Pi= \sum_{I=1}^k \partial_{q_I}\wedge
\partial_{p_I}+ \sum_{i,j=1}^{m-k} \varphi_{ij}(y)\partial_{y_i}\wedge
\partial_{y_j}
\end{equation}
for functions $\varphi_{ij}\colon \RR^{m-k}\rightarrow \RR$.
\end{prop}
\begin{remark}\label{injanc}
The existence of coordinates in which $\Pi$ has the above split
form is guaranteed by Weinstein's Splitting Theorem \cite{We}; the
point in the above proposition is that one can choose the
coordinates $(q,p,y)$ so that $M$ is given by the constrains
$p=0$. \end{remark}

\begin{proof} We adapt the proof of Weinstein's Splitting
Theorem \cite{We} to our setting.
 To simplify the notation we will
often write $P$ in place of $U$ and $M$ in place of $M\cap U$.
We proceed by induction over $k$; for $k=0$ there is nothing to prove, so let $k>0$.

Choose a function $q_1$ on $P$ near $x$ such that $dq_1$ does not
annihilate $\sharp N^*M$. Then $X_{q_1}|_M$  is
transverse to $M$, because there is a $\xi\in N^*M$ with $0\neq
\langle \sharp \xi, dq_1 \rangle=-\langle \xi, X_{q_1} \rangle$.
Choose a hypersurface in $P$ containing $M$ and transverse to
$X_{q_1}|_M$, and determine the function $p_1$ by requiring that
it vanishes on the hypersurface and $dp_1(X_{q_1})=-1$. Since
$[X_{q_1},X_{p_1}]=X_1=0$ the span of $X_{p_1}$ and $X_{q_1}$ is
an integrable distribution giving rise to a foliation of $P$ by
surfaces. This foliation is transverse to $P_1$, which we define as the codimension
two submanifold where $p_1$ and $q_1$ vanish. $M_1:=P_1\cap M$ is a
clean intersection and is a codimension one submanifold of $M$. To
proceed inductively we need

\begin{lemma} \label{ind}
$P_1$ has an induced Poisson structure $\Pi_1$, $M_1\subset P_1$ is a
coisotropic submanifold of codimension $k-1$, and the sharp-map $\sharp_1$ of $P_1$ is injective on
the conormal bundle to $M_1$.
\end{lemma}
\begin{proof}
$P_1$ is cosymplectic because  $\sharp N^*P_1$ is spanned by $X_{q_1}$ and $X_{p_1}$, which are transverse to $P_1$.
Hence it has an induced Poisson structure $\Pi_1$. Recall from section \ref{def} that if $\xi_1\in
T_x^*P_1$ then $\sharp_1\xi_1\in T{P_1}$ is given as follows:
extend $\xi_1$ to a covector $\xi$ of $P$ by asking that it
annihilate $\sharp N_x^*P_1$ and apply $\sharp$ to it. Now in
particular let $x\in M_1$ and $\xi_1$ be
 an element of the conormal bundle of $M_1$ in
$P_1$. We have $T_xM=T_xM_1\oplus\RR X_{p_1}(x)\subset
T_xM_1+\sharp N_x^*P_1$, so $\xi \in N_x^*M$, and since  $M$ is
coisotropic in $P$ we have $\sharp \xi \in T_xM$ . Hence
$\sharp_1\xi_1 \in T_xP_1\cap T_xM=T_xM_1$, which shows the
claimed coisotropicity. The injectivity of $\sharp_1$ on the
conormal bundle follows by the above together with the injectivity
of $\sharp|_{N^*M}$, which holds by Lemma \ref{easylemma}.
\end{proof}

By the induction assumption there are coordinates on $P_1$ so that 
$$ \Pi_1= \sum_{I=2}^k \partial_{q_I}\wedge
\partial_{p_I}+ \sum_{i,j=1}^{m-k} \varphi_{ij}(y)\partial_{y_i}\wedge
\partial_{y_j}$$ and 
$M_1\subset P_1$ is given by the constraints
$p_2=0,\dots,p_k=0$. 
We extend the coordinates on $P_1$ to   the whole of $P$ so that 
they are constant along the surfaces tangent to
$span\{X_{q_1},X_{p_1}\}$. 
We denote collectively by $x_{\alpha}$ the resulting  functions on $P$, which together with $q_1$ and $p_1$ form a coordinate system on $P$.  
  We have $\{x_{\alpha},q_{1}\}=0$ and
$\{x_{\alpha},p_{1}\}=0$, and using the Jacobi identity one sees
that $\{x_{\alpha},x_{\beta}\}$ Poisson commutes with $q_{1}$ and
$p_{1}$, and hence it is a  function of the
$x_{\alpha}$'s only. Further $\{x_{\alpha},x_{\beta}\}|_{P_1}=\{{x_{\alpha}}|_{P_1},{x_{\beta}}|_{P_1}\}_1$ since $x_{\alpha},x_{\beta}$ annihilate $\sharp N^*P_1$.
 Hence formula $\eqref{canform}$
for the Poisson bivector $\Pi$ follows.

To show that $M$ is given by the constraints $p_1=\cdots=p_k=0$
we notice the following.   $p_1$ was chosen to vanish on $M$. The functions 
$p_2,\cdots,p_k$ on $P_1$ were chosen to to vanish on  $M_1$, and since   $TM|_{M_1}=TM_1\oplus \RR X_{p_1}|_{M_1}$ it follows that their extensions
vanish on the whole of $M$.   This concludes the proof of Prop. \ref{can}.
\end{proof}

Using the normal forms derived above we can prove local uniqueness:
\begin{prop}\label{locun}
Suppose we are given a  Dirac manifold $(M,L)$ for which $L\cap
TM$ has constant rank $k$, and let $(P,\Pi)$ be a
Poisson manifold of dimension $\dim M+ k$ in
which $(M,L)$ embeds coisotropically. Then
about each $x\in M$ there is a neighborhood $U\subset P$  
 which is Poisson diffeomorphic
to an open set in the canonical Poisson manifold associated to $(M,L)$ in Prop.  \ref{splitting}.
\end{prop}

\begin{proof}
By Prop. \ref{can}  there are coordinates $\{q_I,p_I,y_i\}$ 
  on $U$ such that locally $M$ is given by  
$p_I=0$ and
$$\Pi= \sum_{I=1}^k \partial_{q_I}\wedge
\partial_{p_I}+ \sum_{i,j=1}^{m-k} \varphi_{ij}(y)\partial_{y_i}\wedge
\partial_{y_j}.$$
We want to apply the construction of Thm. \ref{E} to $(M,L)$. To do so we need to make a 
choice of complement to
$E:=L\cap TM=span\{{\partial_{q_I}}|_M\}$; our choice is $V:=span\{{\partial_{y_i}}|_M\}$.
 Since by
assumption $L$ is the pullback 
  of $graph(\Pi)$ to $M$, $L$ is spanned by sections $({\partial_{q_I}}|_M\oplus 0)$ and 
  $(\sum_j \varphi_{ij}(y){\partial_{y_j}}|_M\oplus {dy_i}|_M)$. 
  Hence the pullback of $L$ to the total space of the vector bundle $\pi:E^* \rightarrow M$ is
  spanned by $({\partial_{q_I}}\oplus 0)$, $({\partial_{p_I}}\oplus 0)$, and 
  $(\sum_j \varphi_{ij}(y){\partial_{y_j}}\oplus {dy_i})$. Next we consider the embedding $E^*\rightarrow
  T^*M$ induced by the splitting $TM=E\oplus V$ and pull back the canonical 2-form
  on $T^*M$. In the coordinates $(q_I,y_i)$ on $M$ the pullback 2-form   is simply
  $\sum_{I=1}^k dp_I\wedge dq_I$ (see eq. (6.7) in \cite{OP}),
  where with $p_I$ we denote the linear coordinates on the fibers of $E^*$ dual to the $q_I$.
  Transforming $\pi^*L$ by this 2-form gives exactly $graph(\Pi)$.
Hence we conclude that, nearby $x\in M$,  the Poisson manifold $(P,\Pi)$  is obtained by the construction
of Thm. \ref{E} (with the above choice of distribution $V$).
\end{proof}

\bibliographystyle{habbrv}
\bibliography{bibmarco}
\end{document}